\documentclass[a4paper,12pt,numbers]{article}
\usepackage{amsmath, amsthm, mathrsfs, amssymb, verbatim, bbm, url}  
\usepackage{amscd}
\usepackage{latexsym,amssymb}
\usepackage{graphicx}
\usepackage{color}
\usepackage{enumerate}
\numberwithin{equation}{section}
    \newcommand{\ga}[0]{\alpha }
\newcommand{\gb}[0]{\beta }

\newcommand{\gd}[0]{\delta }
\newcommand{\gf}[0]{\varphi }

      \newcommand{\la}{\left\langle}
      \newcommand{\xii}{\mathcal{X}_{g,n,d,(\mu_1,\ldots ,\mu_n)}}
      \newcommand{\xiii}{\mathcal{X}_{g,n+1,d,(\mu_1,\ldots ,\mu_n,0)}}

        \newcommand{\tzz}{\widetilde{\mathcal{Z}}}
     \newcommand{\ops}{\overline{\psi}}
    \newcommand{\ra}{\right\rangle}
     \newcommand{\ozz}{\mathcal{O}_\mathcal{Z}}
         \newcommand{\ozzr}{\mathcal{O}_{\mathcal{Z}_\mu}}
         \newcommand{\di}{\mathcal{D}_{j,(\mu_1,\ldots ,\mu_n)}}
    \newcommand{\ix}{\mathcal{X}}
      
     \newcommand{\oddi}{\mathcal{O}_{\mathcal{D}_j}}
     \newcommand{\ic}{\mathcal{U}}

\newcommand{\BB}{\mathcal{B}_{g,n,d}}
\newcommand{\CC}{\mathcal{C}_{g,n,d}}
\newcommand{\AAA}{\mathcal{A}_{g,n,d}}
\newtheorem{loopg}{Theorem}[section]
\newtheorem{dil}[loopg]{Theorem}
\newtheorem{polaris}[loopg]{Theorem}
\newtheorem{mot11}[loopg]{Remark}
\newtheorem{tthsian}[loopg]{Remark}
\newtheorem{rem35}[loopg]{Remark}
\newtheorem{rem37}[loopg]{Remark}
\newtheorem{ttpol}[loopg]{Remark}
\newtheorem{gworb}{Definition}[section]
\newtheorem{gworb2}[gworb]{Definition}
\newtheorem{gwex}[gworb]{Example}
\newtheorem{gworb5}[gworb]{Remark}
\newtheorem{useless}[gworb]{Remark}
\newtheorem{gworb4}[gworb]{Theorem}
\newtheorem{rem3}[gworb]{Remark}
\newtheorem{altarem}[gworb]{Remark}
\newtheorem{teta}{Proposition}[section]
\newtheorem{pushpull}[teta]{Lemma}
\newtheorem{tetaaa}[teta]{Proposition}
\newtheorem{tetabb}[teta]{Proposition}
\newtheorem{exbet}[teta]{Example}
\newtheorem{chernch}[teta]{Lemma}
\newtheorem{rem11}[teta]{Remark}
\newtheorem{nodalclass4}[teta]{Lemma}
\newtheorem{tetacc}[teta]{Proposition}
\newtheorem{nodalclass2}[teta]{Lemma}
\newtheorem{coord1}{Remark}[section]
\newtheorem{coord2}[coord1]{Example}
\newtheorem{defcoh}{Definition}[section]
\newtheorem{indexx0}[defcoh]{Theorem}
\newtheorem{tangvirtual}[defcoh]{Proposition}
\newtheorem{indexx}[defcoh]{Corollary}
\newtheorem{typeaa}{Corollary}[section]
\newtheorem{typebb}[typeaa]{Corollary}
\newtheorem{typecc}[typeaa]{Corollary}
\newtheorem{trace}{Definition}[section]
\newtheorem{cher}[trace]{Definition}
\newtheorem{todd}[trace]{Definition}
\newtheorem{Grr}[trace]{Theorem}

\begin{document}

\begin{center}
\Large{ Twisted orbifold Gromov-Witten invariants}\footnote{ Mathematics Subject Classification 14N35} \vspace{1 cm}

\large{Valentin Tonita} 
\end{center}
\begin{abstract}
Let $\ix$ be a smooth proper Deligne-Mumford stack over $\mathbb{C}$. One can define twisted orbifold Gromov-Witten invariants of $\ix$ by considering multiplicative invertible characteristic classes of various bundles on the moduli spaces of stable maps $\ix_{g,n,d}$, cupping them with evaluation and cotangent line classes and then integrating against the virtual fundamental class. These are more general than the twisted invariants introduced in \cite{tseng}. We express the generating series of the twisted invariants in terms of the generating series of the untwisted ones. We derive the corollaries which are used in the paper \cite{gito} on the quantum K-theory of a complex compact manifold $X$.
 \end{abstract}
 \tableofcontents
\section{Introduction and statement of results}

 Twisted Gromov-Witten invariants have been introduced in \cite{coatesgiv} for manifold target spaces and extended by \cite{tseng} to the case of 
 orbifolds. The original motivation was to express Gromov-Witten invariants of complete intersections (the ``twisted'' ones) in terms of the GW invariants of the ambient space (the untwisted ones). In addition they were used in \cite{Tom} to express Gromov-Witten invariants with values in cobordism in terms of cohomological Gromov-Witten invariants. 
 
  Our results incorporate and generalize all of the above: we consider three types of twisting classes. These are multiplicative cohomological classes of bundles of the form $\pi_*E$, where $\pi$ is the universal family of the moduli space of stable maps to an orbifold $\ix$. The main tool in the computations is the Grothendieck-Riemann-Roch theorem for stacks of \cite{toen}, applied to the morphism $\pi$: this gives differential equations satisfied by the generating functions of the twisted GW invariants. To the genus $0$ Gromov-Witten potential of an orbifold $\ix$ one can associate an overruled Lagrangian cone in a symplectic space $\mathcal{H}$ - as explained in Section 2. Solving the differential equations for each type of twisting has an interpretation in terms of the geometry of the cone: change its position  by a symplectic transformation, translation of the origin and a change of polarization of $\mathcal{H}$. Our motivation comes from studying the quantum K-theory of a manifold $X$ (see \cite{gito}), more precisely trying to express Euler characteristics on the (virtual) orbifolds $X_{0,n,d}$ in terms of cohomological Gromov-Witten invariants. However they have other applications - for instance recovering the work of \cite{Tom} on quantum extraordinary cohomology.
  
  In \cite{tel}, Teleman studies a group action on 2 dimensional quantum field theories. Our results match his, if the field theories come from Gromov-Witten theory.
   
 Let $\ix$ be a compact orbifold. Moduli spaces of orbimaps to orbifolds have been constructed by \cite{chru} in the setup of symplectic orbifolds and by \cite{abgrvi2} in the context of Deligne-Mumford stacks. Informally, the domain curve is allowed to have nontrivial orbifold structure at the marked points and nodes. We denote the moduli spaces of degree $d$ maps of genus $g$  with $n$ marked points by $\ix_{g,n,d}$.
 
 Just like in the case of manifold target spaces, there are evaluation maps $\overline{ev}_i$ at the marked points. Although it is clear how these maps are defined on geometric points, it turns out that to have well-defined morphisms of Deligne-Mumford stacks the target of the evaluation maps is the rigidified inertia stack of $\ix$. We first define a related object, the inertia stack $I\ix$, as follows: around any point $x\in \ix$ there is a local chart $(\widetilde{U}_x, G_x)$ such that locally $\ix$ is represented as the quotient of $\widetilde{U}_x$ by $G_x$. Consider the set of conjugacy classes $(1)=(h_x^1)$, $(h_x^2)$, $\ldots$, $(h_x^{n_x})$ in $G_x$. Define:
      \begin{align*}
      I\ix:=\{(x,(h_x^i)) \quad \vert \quad i=1,2,\ldots, n_x \}.
       \end{align*}
   Pick an element $h_x^i$ in each conjugacy class. Then a local chart on $I\ix$ is given by:
    \begin{align*}
       \coprod_{i=1}^{n_x} \widetilde{U}_x^{(h_x^i)}/ Z_{G_x}(h_x^i),
     \end{align*}
     where $Z_{G_x}(h_x^i)$ is the centralizer of $h_x^i$ in $G_x$.

   The rigidified inertia stack, which we denote $\overline{I\ix}$, is defined by taking the quotient at $(x,(g))$ of the automorphism group by the cyclic subgroup generated by $g$. So, whereas a local chart at $(x,(g))$ on $I\ix$ is given by $\widetilde{U_g}/Z_{G_x}(g)$, on $\overline{I\ix}$ a local chart is $\widetilde{U_g}/[Z_{G_x}(g)/\langle g\rangle ]$. It is in general disconnected, even if $\ix$ is connected. We write $\overline{I\ix}:=\coprod_\mu \overline{\ix_\mu}$. The distinguished component corresponding to the identity is a copy of $\ix$ and throughout we will label it $\ix_0$ to distinguish it from other components of $I\ix$. We denote by $\iota : I\ix\to I\ix$ the involution which maps $(x,(g))$ to $(x,(g^{-1}))$. It descends to an involution on $\overline{I\ix}$, which we also denote $\iota$. We write $\ix_{\mu^I}: = \iota(\ix_\mu)$. There is a natural map $q:I\ix\to \ix$.
   
    The orbifold Poincar\'e pairing on $I\ix$ is defined for $a \in H^*(\ix_\mu, \mathbb{C})$ , $b \in H^*(\ix_{\mu^I}, \mathbb{C})$ as:
      \begin{align*}
      (a,b)_{orb}: = \int_{\ix_\mu} a\cup \iota^* b .
      \end{align*}
   
   $I\ix$ and $\overline{I\ix}$ have the same geometric points (coarse spaces), hence we can identify the rings $H^*(I\ix,\mathbb{C})$ and $H^*(\overline{I\ix},\mathbb{C})$. This allows us to pretend the cohomological pullbacks by the maps $\overline{ev}_i$ have domain $H^*(I\ix ,\mathbb{C})$. 
 We can use the maps $\overline{ev}_i$ to decompose $\ix_{g,n,d}$ as a union of closed and open substacks:
  \begin{align*}
  \xii : = \ix_{g,n,d}\cap (\overline{ev}_1)^{-1}(\overline{\ix}_{\mu_1})\cap\ldots \cap (\overline{ev}_n)^{-1}(\overline{\ix}_{\mu_n})
   \end{align*}
  
    For each $i$ we denote by $\ops_i = c_1 (\overline{L}_i)$, where the line bundle  $\overline{L}_i$ has fiber over each point $(\mathcal{C},x_1,\ldots,x_n,f)$  the cotangent line to the coarse curve $C$ at $x_i$.  
  
    We denote the universal family by $\pi : \mathcal{U}_{g,n,d}\to \ix_{g,n,d}$. $\mathcal{U}_{g,n,d}$ can be identified with $\cup_{(\mu_1,\ldots ,\mu_n)}\xiii$. Since the extra marked point on the universal family has trivial orbifold structure the map $\overline{ev}_{n+1}$ lands in $\ix_0$. We will write $ev_{n+1}$ throughout.
  The moduli spaces $\ix_{g,n,d}$ have  perfect obstruction theory and are equipped with virtual fundamental classes $[\ix_{g,n,d}]\in H_*(\ix_{g,n,d},\mathbb{Q})$.
   Orbifold Gromov-Witten invariants are obtained by integrating $\ops_i$ and evaluation classes on these cycles. We use correlator notation:
    \begin{align*}
  \la a_1\ops^{k_1},\ldots ,a_n \ops^{k_n}\ra_{g,n,d}: = \int_{[\ix_{g,n,d}]}\prod^n_{i=1}ev_i^*a_i \ops_i^{k_i} . 
     \end{align*}    
    
   Their generating series are functions on a suitable infinite dimensional vector space $\mathcal{H}_+$, which we describe below.
 Let  $\Lambda: = \mathbb{C}[[Q]]$ be the Novikov ring which is a completion of the semigroup ring of degrees of holomorphic curves in $\ix$ and let:  
   \begin{align*}
 \mathcal{H}:=  H^*(I\ix,\Lambda)((z)) .
    \end{align*}
 We equip $\mathcal{H}$ with the symplectic form:
   \begin{align*}
  \Omega(\mathbf{f},\mathbf{g}):= \oint_{z=0}\left(\mathbf{f}(z), \mathbf{g}(-z)\right)_{orb}dz  .
   \end{align*}    
  Consider the following polarization of $\mathcal{H}$:
   \begin{align*}
  \mathcal{H}_+:= H^*(I\ix,\mathbb{C})[[z]] \quad and \quad \mathcal{H}_-:= z^{-1} H^*(I\ix,\mathbb{C})[z^{-1}].
   \end{align*}  
 Let $\mathbf{t}(z)\in \mathcal{H}_+$. The genus $g$  descendant potential, respectively the total descendant potential are defined as:
   \begin{align*}
   &\mathcal{F}_\ix^g (\mathbf{t}(z))=\sum_{d,n}\frac{Q^d}{n!}\la \mathbf{t}(\ops),\ldots , \mathbf{t}(\ops)\ra_{g,n,d},\\
   &\mathcal{D}_\ix(\mathbf{t})=exp\left( \sum_{g\geq 0} \hbar^{g-1} \mathcal{F}^g(\mathbf{t})\right).
    \end{align*}   
  Then $\mathcal{D}_\ix$ is a well defined formal function on $\mathcal{H}_+$ taking values in $\Lambda\otimes \mathbb{C}[[\hbar,\hbar^{-1}]]$.  
Also it is well-known that the differential of the genus $0$ potential gives rise to a cone $\mathcal{L}^H\subset \mathcal{H}$ with nice geometric properties (see Theorem $\ref{conethm}$).  
 
  ``Twisted Gromov-Witten invariants'' are obtained from the usual ones by systematically inserting in the correlators multiplicative classes of certain bundles. For a vector bundle $E$, a general multiplicative class is of the form
   \begin{align*}
  \mathcal{A}(E) =exp\left(\sum_{k\geq 0} s_k ch_k E  \right) .
     \end{align*} 
    We want to consider three types of twistings, each by several possibly different multiplicative characteristic classes:
    \begin{itemize}
   \item twistings by a finite number of multiplicative classes $\mathcal{A}_\alpha (\pi_*(ev_{n+1}^* E))$, where $E\in K^0(\ix)$.
   \item twistings by classes $\mathcal{B}_\beta$ (kappa classes) of the form:
     \begin{align*}
     \BB =  \prod_{\beta=1}^{i_B} \mathcal{B}_\beta \left( \pi_*(f_{\beta}(L_{n+1}^{-1})-f_{\beta}(1)) \right) ,
     \end{align*}
    where $L_{n+1}$ is the cotangent line bundle at the extra marked point on the universal curve, $f_\gb$ are polynomials with coefficients in $ev^*_{n+1} K^0(\ix)$ and $1$ is the trivial line bundle.
    \item twistings by nodal classes $\mathcal{C}_\gd$ of the form:
      \begin{align*}
      \CC = \prod_\mu \prod_{\gd=1}^{i_{\mu}} \mathcal{C}^\mu_\gd \left( \pi_* (ev_{n+1}^* F_{\gd\mu} \otimes i_{\mu*}\mathcal{O}_{\mathcal{Z}_\mu})\right),
      \end{align*}
    where $F_{\gd\mu} \in K^0(\ix)$. See Section $2$ for the precise definition of $\mathcal{Z}_\mu$ - roughly speaking it parametrizes nodes with fixed orbifold type; we denote by $i_\mu$ the corresponding inclusion $\mathcal{Z}_\mu\to \mathcal{U}_{g,n,d}$. Hence we allow different types of twistings localized near the loci $\mathcal{Z}_\mu$. 
     \end{itemize}
   We will refer to these as type $\mathcal{A,B,C}$ twistings respectively.  So a twisted GW invariant will be an integral of the form
   \begin{align*}
   \int_{[\ix_{g,n,d}]}\prod^n_{i=1}ev_i^*a_i \ops_i^{k_i}\mathcal{A}(\cdot)\mathcal{B}(\cdot)\mathcal{C}(\cdot).
   \end{align*}
  These can be packed in generating series - the twisted potentials $\mathcal{F}^g_{\mathcal{A,B,C}}$, $\mathcal{D}_{\mathcal{A,B,C}}$, which we can regard as functions on the same space $\mathcal{H}_+$. We postpone the precise definitions to Section 2. We will write $\mathcal{D}_{\mathcal{A,B}}, \mathcal{L}_\mathcal{A}$  etc. for objects associated to twisted GW invariants of the types specified in notation. 

The main theorems of the paper describe how the twistings change the potentials and the corresponding Lagrangian cones $\mathcal{L}_{\mathcal{A,B,C}}$ (which we define in Section 2).        
 \begin{loopg} \label{loop11} 
The cone $\mathcal{L}_\mathcal{A}$ is obtained from $\mathcal{L}^H$ after rotation by a symplectic transformation
    \begin{align*}
    \mathcal{L}_\mathcal{A} = \left(\prod_{\alpha} \Delta_\alpha\right) \mathcal{L}^H.
    \end{align*} 
  We will write explicit formulas for each $\Delta_\alpha$ in Remark $\ref{1090}$.   
     \end{loopg} 
 
Let now $\mathbf{L}_z$ be a line bundle with first Chern class $z$.
   
 \begin{dil} \label{dilat}
The twisting by the classes $\BB$  has the same effect as a translation on the Fock space:
  \begin{align}
\mathcal{D}_{\mathcal{A,B,C}}(\mathbf{t}) =\mathcal{D}_{\mathcal{A,C}}\left( \mathbf{t} +z - z\prod_{i=1}^{i_B}\mathcal{B}_\beta\left( -\frac{f_\beta(\mathbf{L}_{z}^{-1})-f_\beta(1)}{\mathbf{L}_{z}-1}\right)\right) \cdot K_{\mathcal{B}} , \label{bbpp}
  \end{align}
    \end{dil}
  where $K_\mathcal{B}$ is a constant discussed in the proof.
  
 A related result for manifold target spaces is in the paper \cite{kakim}.

\begin{polaris} \label{pola11}
 The potential $\mathcal{D}_{\mathcal{A,B,C}}$ satisfies the differential equation
 \begin{align*}
  \mathcal{D}_{\mathcal{A,B,C}} = exp\left(\frac{\hbar}{2}\sum_{a,b,\ga,\gb,\mu}A^\mu_{a,\alpha;b;\beta}\partial_{a}^{\alpha,\mu}\partial_b^{\beta,\mu^I} \right)\mathcal{D}_{\mathcal{A,B}},
\end{align*} 
where the coefficients $A^\mu_{a,\alpha;b;\beta}$ are defined by formula $\ref{polcoeff}$ in Section 4. This is equivalent to considering the potential  $\mathcal{D}_{\mathcal{A,B}}$ as a generating function with respect to a new polarization  $\mathcal{H}=\mathcal{H}_+\oplus \mathcal{H}_{-,\mathcal{C}}$. We give a precise linear transformation of Darboux coordinates on $\mathcal{H}$ in formula ($\ref{coordinates}$).
      \end{polaris}
 
 A few remarks are in order at this point:
 \begin{mot11}
 \em{The study of the K-theoretic GW invariants of a manifold $X$ in \cite{gito} leads naturally to considering these twisted GW invariants. Briefly put, to compute K-theoretic GW invariants of $X$ in terms of cohomological ones one needs to consider cohomological integrals twisted by certain Todd-like classes (see Section 6) of the (virtual) tangent bundle of $X_{0,n,d}$. Proposition $\ref{tggv}$ expresses this tangent bundle  as a sum of three contributions - one of each type. }
 \end{mot11}

\begin{tthsian}\label{1090}
\em{Theorem $\ref{loop11}$ is a rather straight-forward generalization of the main theorem in \cite{tseng}, the only difference being that we consider more than one class $\mathcal{A}_\alpha$. If the twisting data $\mathcal{A}$ is given by the multiplicative class $ \mathcal{A}(\cdot) = exp (\sum s_k ch_k (\cdot )) $ and by $E\in K^0(\ix )$ then the symplectic transformation $\Delta$ is defined as}
\end{tthsian}
 
  \begin{align*}
 \Delta:= exp\left(\sum_{k\geq 0} s_k\left( \sum_{m\geq 0} \frac{(A_m)_{k+1-m}z^{m-1}}{m!}+ \frac{ch_k (E^{(0)})}{2}\right)\right),
   \end{align*} 
  where by $(A_m)_{j}$ we mean the degree $j$ part of operators of ordinary multiplication by certain elements  $A_m\in H^*(I\ix)$. To define $A_m$ we introduce more notation:  let $r_\mu$ be the order of each element in the conjugacy class which is labeled by $\ix_\mu$. The restriction of the bundle $E$ to $\ix_{\mu}$ decomposes into characters : let $E^{(l)}_{\mu}$ be the subbundle on which every element of the conjugacy class acts with eigen value $e^{2\pi i l/r_\mu}$. The Bernoulli polynomials $B_m(x)$ are defined by 
 \begin{align*}
 \frac{te^{tx}}{e^t-1}=\sum_{m\geq 0}\frac{B_m(x)t^m}{m!}.
 \end{align*}
  Then
   \begin{align*}
  (A_m)_{\vert \ix_{\mu} }:= \sum^{l=r_\mu -1}_{l=0}B_m(\frac{l}{r_\mu}) ch(E^{(l)}_{\mu}) .
   \end{align*} 
The symplectic operator in Theorem $\ref{loop11}$ is just the product of Tseng's operators $\Delta_\alpha$ associated to each $\mathcal{A}_\alpha$. 
 
    \begin{rem35}
   {\em  The decomposition:
   \begin{align*}
  H^*(I\ix,\mathbb{C})((z^{-1}))=\oplus H^*(\ix_{\mu},\mathbb{C}) ((z^{-1}))
   \end{align*} 
   is preserved by the action of this loop group element. $A_m$ acts by cup product multiplication on each $H^*(\ix_{\mu})$}.
   \end{rem35}     
  \begin{rem37} \label{rem333}
  {\em Theorem $\ref{loop11}$ can be extended to a statement about the total descendant potential using the quantization formalism of \cite{gi}.  It reads:}
   \begin{align*}
  \mathcal{D}_\mathcal{A}(\mathbf{q})  \approx \prod_{\alpha} \widehat{\Delta}_\alpha \mathcal{D}_\ix(\mathbf{q} ) ,
   \end{align*} 
{\em where $\widehat{\Delta}$ denotes the quantization of the operator $\Delta$ and the symbol $\approx$ means the two sides are equal up to a (precisely determined) scalar factor. }
  \end{rem37}
  
\begin{ttpol}
\em{Another way to obtain a basis for the new space $\mathcal{H}_{-,\mathcal{C}}$ of the new polarization from the Theorem $\ref{pola11}$ is the following: for each $\mu$ let the series $u_\mu(z)$ be defined by}
 \begin{align*}
 \frac{z}{u_\mu(z)}= \prod^{i_\mu}_{\gd=1}\mathcal{C}^\mu_\gd\left((q^*F_{\gd\mu})^{(0)}_\mu \otimes(-\mathbf{L}_{-z}) \right). 
 \end{align*}
Moreover define Laurent series $v_{k,\mu}$, $k=0,1,2,\ldots$ by:
 \begin{align*}
  \frac{1}{u_\mu(-x-y)} =\sum_{k\geq 0}(u_\mu(x))^k v_{k,\mu}(u(y)) \quad,   
 \end{align*} 
    where we expand the left hand side in the region where $\vert x\vert<\vert y\vert$ .   
 Then  $\mathcal{H}_{-,\mathcal{C}}=\oplus_\mu \mathcal{H}^\mu_{-,\mathcal{C}}$  and each $\mathcal{H}^\mu_{-,\mathcal{C}}$ is spanned by $\{\varphi_{\ga,\mu}v_{k,\mu}(u(z))\}$ where $\{\varphi_{\ga,\mu}\}$ runs over a basis of $H^*(\ix_\mu,\mathbb{C})$ and $k$ runs from $0$ to $\infty$.    
 \end{ttpol} 
 
   The rest of the paper is structured as follows. Section 2 is used to introduce the  main  objects of study: the moduli spaces $\ix_{g,n,d}$ and the Gromov-Witten theory of $\ix$, the symplectic space $\mathcal{H}$, the (twisted and untwisted) Gromov-Witten potentials. Section 3 contains the technical results which are the core of the computations - mainly how the twisting cohomological classes pullback on the universal family and the locus of nodes. We are now ready to prove the Theorems $\ref{loop11}$, $\ref{dilat}$ and $\ref{pola11}$ - which we do in Section 4. In Section 5 we use the results to give a concise proof of the fake quantum Hirzebruch-Riemann-Roch theorem: this was done in \cite{Tom} by a very long calculation. In Section 6 we extract the corollaries which are used in the paper \cite{gito} on quantum K-theory.  
Finally, in the appendix we state To\"en's Grothendieck-Riemann Roch theorem for stacks, which applied to the universal family is the starting point in the computation.  

$ \mathbf{Acknowledgments}$. I would like to thank Alexander Givental for suggesting the problem as a tool for the work in \cite{gito} and to Tom Coates and Hsian-Hua Tseng for useful discussions.
 
    \section{Orbifold Gromov-Witten theory} 
    Throughout this paper, $\ix$ will be a proper smooth Deligne-Mumford stack over $\mathbb{C}$ with projective coarse moduli space.
  
     We now recall the definitions of orbicurve and of orbifold stable maps of \cite{chru} and \cite{abgrvi2}. The idea to extend the definition of a stable map to an orbifold target space is quite natural. One then notices that in order to obtain compact moduli spaces parametrizing these objects one has to allow orbifold structure on the domain curve at the nodes and marked points (see e.g. \cite{ab1}).   
   \begin{gworb}
 {\em A nodal $n$-pointed orbicurve $(\mathcal{C}, x_1 , x_2, \ldots , x_n)$ is a nodal marked complex curve such that } 
   \end{gworb}   
  \begin{itemize}
  \item $\mathcal{C}$ has trivial orbifold structure on the complement of the marked points and nodes.
  \item  Locally near a marked point,  $\mathcal{C}$ is isomorphic to $[$Spec $\mathbb{C}[z]/\mathbb{Z}_r]$, for some $r$, and the generator of $\mathbb{Z}_r$ acts by $z\mapsto \zeta z$, $\zeta^r=1$. 
  \item  Locally near a node,  $\mathcal{C}$ is isomorphic to  $[$Spec $\left(\mathbb{C}[z,w]/(zw)\right) /\mathbb{Z}_r]$,  and the generator of $\mathbb{Z}_r$ acts by $z\mapsto \zeta z$, $w\mapsto \zeta^{-1}w$. We call this action {\em balanced} at the node.
  \end{itemize} 
  We now define twisted stable maps:
  \begin{gworb2} \label{def111}
 {\em An  $n$-pointed, genus $g$, degree $d$ orbifold stable map is a representable morphism $f: \mathcal{C}\to \ix$ , whose domain is an $n$-pointed genus $g$ orbicurve $\mathcal{C}$ such that the induced morphism of the coarse moduli spaces $C\to X$ is a stable map of degree $d$}.  
  \end{gworb2}
 
 We denote the moduli space parametrizing $n$-pointed, genus $g$, degree $d$ orbifold stable maps by $\ix_{g,n,d}$. It is proved in \cite{abvi} that $\ix_{g,n,d}$ is a proper Deligne-Mumford stack. Just like the case of stable maps to manifolds, there are evaluation maps at the marked points, but these land naturally in the {\em rigidified} inertia orbifold of $\ix$, which we denote $\overline{I\ix}$.
 
  \begin{gwex}
  \emph{If $\ix$ is a global quotient $Y/G$ then the strata of $I\ix$ are $Y^g/C_G(g)$ and of $\overline{I\ix}$ are $\overline{\ix}_{(g)} : = Y^g/\overline{C_G(g)}$, where $\overline{C_G(g)}=C_G(g)/\langle g\rangle$ for each conjugacy class $(g)\subset G$. } 
   \end{gwex}  
     
  See \cite{abgrvi} and \cite{abgrvi2} for the definition of $\overline{I\ix}$ in the category of stacks.
   
We decompose $\ix_{g,n,d}$ according to the target of the evaluation maps:
  \begin{align*}
  \xii : = \ix_{g,n,d}\cap (\overline{ev}_1)^{-1}(\overline{\ix}_{\mu_1})\cap\ldots \cap (\overline{ev}_n)^{-1}(\overline{\ix}_{\mu_n}).
   \end{align*}
 Since we work with cohomology with complex coefficients we consider the cohomological pullbacks by the maps $\overline{ev}_i$ having domain $H^*(I\ix ,\mathbb{C})$.  $I\ix$ and $\overline{I\ix}$ have the same coarse spaces, which implies that both spaces have the same cohomology rings with rational coefficients. In fact there is a map $\Pi:I\ix\to \overline{I\ix}$, which maps a point $(x,(g))$ to $(x,(\overline{g}))$. If $r_i$ is the order of the automorphism group of $x_i$, then define:  
    \begin{align*}
   & ev^*_i : H^*(I\ix,\mathbb{C}) \to H^*(\ix_{g,n,d},\mathbb{C}),  \\
   & a\mapsto r_i^{-1} (\overline{ev}_i)^* (\Pi_* a ) .
    \end{align*}   
    
 Notice that if a marked point $x_i$ has trivial orbifold structure, $\overline{ev}_{i}$ lands in  the distinguished component $\ix_0$ of $\overline{I\ix}$. The universal family can be therefore identified with the diagram:
 \begin{displaymath}
     \begin{CD}
    \ic_{g,n,d}:= \cup_{(\mu_1,\ldots , \mu_n)}\xiii @>ev_{n+1}>> \ix \\
      @V \pi VV        \\
       \ix_{g,n,d}               \qquad .
     \end{CD}
     \end{displaymath} 
  
 In the universal family $\ic_{g,n,d}$ lies the divisor of the $i$-th marked point $\mathcal{D}_i$: its points parametrize maps whose domain has a distinguished node separating two orbicurves $\mathcal{C}_0$ and $\mathcal{C}_1$. $\mathcal{C}_1$ has genus $0$ and carries only three special points: the node, the $i$-th marked point and the $(n+1)$-st marked point and is mapped with degree $0$ to $\ix$. We write:
  \begin{align*}
   \mathcal{D}_{i,(\mu_1,\ldots ,\mu_n)}:= \mathcal{D}_i \cap \xiii  .   
   \end{align*}
  We denote by $\sigma_i$ the corresponding inclusions.
 
  Let $\mathcal{Z}$ be the locus of nodes in the universal family. It has codimension two in $\ic_{g,n,d}$.  Denote by $p:\tzz\rightarrow \mathcal{Z}$ the double cover  over $\mathcal{Z}$ given by a choice of $+,-$ at the node. For the inclusion of a stratum:
\begin{align*}
   \ix_{g_1,n_1+1,d_1}\times_{\overline{I\ix}} \ix_{0,3,0}\times_{\overline{I\ix}}\ix_{g_2,n_2+1,d_2}\rightarrow\mathcal{Z} \hookrightarrow  \ix_{g,n+1 ,d}
     \end{align*}
   we  will denote by $p_i$ ($i=1,2$) the projections:
     \begin{align*}
  p_i : \ix_{g_1,n_1+1,d_1}\times_{\overline{I\ix}} \ix_{0,3,0}\times_{\overline{I\ix}}\ix_{g_2,n_2+1,d_2} \rightarrow  \ix_{g_i,n_i+1,d_i}.
     \end{align*}
  
  We denote $\mathcal{Z}^{irr}, \mathcal{Z}^{red} $ the loci of nonseparating nodes, respectively separating nodes and  $i^{irr}, i^{red}$ for the inclusion maps. Moreover we will need to keep track of the orbifold structure of the node. We denote by $\mathcal{Z}_\mu$ the locus of nodes where the evaluation map at one branch lands in $\overline{\ix }_\mu$ and by $i_\mu$ the corresponding inclusions. 
    
   The moduli spaces $\ix_{g,n,d}$ have perfect obstruction theory (see \cite{abgrvi2}). According to \cite{behfan} this yields virtual fundamental classes:
   \begin{align*} 
    [\ix_{g,n,d}]\in H_*(\ix_{g,n,d} ,\mathbb{Q}).
    \end{align*}
    
   We define $\ops_i$ to be the first Chern classes of line bundles whose fibers over each point $(\mathcal{C},x_1,\ldots ,x_n , f)$  are  the cotangent spaces at $x_i$ to the \emph{coarse curve} $C$. GW invariants are obtained by intersecting $\ops$ and evaluation classes against the virtual fundamental class. We write:

    \begin{align*}
   \langle a_1\ops^{k_1},\ldots , a_n\ops^{k_n} \rangle_{g,n,d}: = \int_{[\ix_{g,n,d}]} \prod_{i=1}^n ev_i^*(a_i)\ops^{k_i}_i .  
    \end{align*}
  \begin{gworb5}
  \emph{ The moduli spaces $\ix_{g,n,d}$ and the evaluation maps, differ from those considered in \cite{tseng}. However the Gromov-Witten invariants agree, since integration in \cite{tseng} is done over a weighted virtual fundamental class. } 
  \end{gworb5} 
  
 Let $\mathbb{C}[[Q]]$ be the Novikov ring which is the formal power series completion of the semigroup ring of degrees of holomorphic curves in $X$. For more on Novikov rings see \cite{mcduff}. We define the ground ring $\Lambda:=\mathbb{C}[[Q]]$ and:   
   \begin{align*}
 \mathcal{H}:=  H^*(I\ix,\Lambda)((z)).
    \end{align*}
   We endow $\mathcal{H}$ with the symplectic form:
   \begin{align*}
  \Omega(\mathbf{f},\mathbf{g}):= \oint_{z=0}\left(\mathbf{f}(z), \mathbf{g}(-z)\right)_{orb}dz  . 
   \end{align*}    
   The following polarization of $\mathcal{H}$:
   \begin{align*}
  \mathcal{H}_+:= H^*(I\ix,\Lambda)[[z]]; \quad  \quad \mathcal{H}_-:= z^{-1} H^*(I\ix,\Lambda)[z^{-1}].
   \end{align*} 
identifies $\mathcal{H}$ with $T^*\mathcal{H}_+$.
\begin{useless}
{\em This choice of polarization is different from the one in most places in literature. The reason is because in applying these results to quantum $K$-theory we need that $e^z\in \mathcal{H}_+$. See \cite{gito} for details.}
\end{useless}
 Let $\{\varphi_\alpha\}$ and $\{\varphi^\beta\}$ be dual bases in $H^*(I\ix,\Lambda)$. We introduce Darboux coordinates $\{p_a^\alpha , q_b^\beta\}$ on $\mathcal{H}$ and we write:
  \begin{align*}
  &\mathbf{p}(z) =\sum_{a,\alpha}p_a^{\alpha}\gf_\alpha (-z)^{-a-1} \in \mathcal{H}_-  ,\\
  & \mathbf{q}(z) =\sum_{b,\gb}q_b^{\gb}\gf^\gb z^{b} \in \mathcal{H}_+ .
  \end{align*}  
  We equip  $\mathcal{H}$ with the $Q$-adic topology. Let:
  \begin{align*}
  \mathbf{t}(z):=t_0 + t_1 z +\dots \in H^*(I\ix,\Lambda)[[z]] .
  \end{align*}  
  
   Then the genus $g$, respectively total potential are defined to be:
   \begin{align*}
   &\mathcal{F}^g (\mathbf{t}(z))=\sum_{d,n}\frac{Q^d}{n!}\la \mathbf{t}(\ops),\ldots , \mathbf{t}(\ops)\ra_{g,n,d} ,\\
   &\mathcal{D}(\mathbf{t}(z))=exp\left( \sum_{g\geq 0} \hbar^{g-1} \mathcal{F}^g(\mathbf{t}(z))\right).
    \end{align*}    

 For $\mathbf{t}(z)\in \mathcal{H}_+$ we call the translation $\mathbf{q}(z):=\mathbf{t}(z)-\mathbf{1}z$ the {\em dilaton shift}. We regard the total descendant potential as a formal function on $\mathcal{H}_+$ in a neghborhood of $-\mathbf{1}z$  taking values in $\mathbb{C}[[Q,\hbar, \hbar^{-1}]]$. 
  
  The graph of the differential of $\mathcal{F}^0$ defines a formal germ of a Lagrangian submanifold of $\mathcal{H}$:
    \begin{align*}
   \mathcal{L}^H:=\{(\mathbf{p},\mathbf{q}), \mathbf{p}=d_{\mathbf{q}}\mathcal{F}^0\}\in \mathcal{H}. 
    \end{align*}
   \begin{gworb4} \label{conethm}
 \emph{(\cite{Given1})}$\mathcal{L}^H$ is (the formal germ of) a Lagrangian cone with vertex at the origin such that each tangent space $T$ is tangent to $\mathcal{L}^H$ exactly along $zT$.   
   \end{gworb4}  
    
   The class of cones satisfying properties of Theorem $\ref{conethm}$ is preserved under the action of symplectic transformations on $\mathcal{H}$ which commute with multiplication by $z$. We call these symplectomorphisms \emph{loop group elements}. They are matrix valued Laurent series in $z$:
     \begin{align*}
   S(z)=\sum_{i\in\mathbb{Z}} S_i z^i  ,
     \end{align*}
     where $S_i \in End\left(H^*(I\ix)\otimes \Lambda\right)$. Being a symplectomorphism amounts to:
      \begin{align*}
    S(z)S^*(-z) = I  ,
     \end{align*} 
   where $I$ is the identity matrix and $S^*$ is the adjoint of $S$. Differentiating the relation above at the identity, we see that infinitesimal loop group elements $R$ satisfy:      
    \begin{align*}
   R(z) + R^*(-z) = 0 . 
    \end{align*}    
   
    We now introduce twisted Gromov-Witten invariants. For a bundle $E$ we will denote by $\mathcal{A}(E)$, $\mathcal{B}(E)$, $\mathcal{C}(E)$ general multiplicative classes of $E$. These are of the form:
      \begin{align*}
      exp\left(\sum_{k\geq 0} s_k ch_k(E)\right)  .  
      \end{align*}
       We then define the classes $\AAA, \BB, \CC \in H^*(\ix_{g,n,d})$ as products of possibly \emph{different} multiplicative classes of bundles:
    \begin{align*}
  & \AAA  = \prod_{\alpha=1}^{i_A} \mathcal{A}_\alpha (\pi_*(ev^* E_\alpha) ) ,\\
   & \BB =  \prod_{\beta=1}^{i_B} \mathcal{B}_\beta \left( \pi_*(f_{\beta}(L_{n+1}^{-1})-f_{\beta}(1)) \right) , \\
   & \CC = \prod_\mu\prod_{\gd=1}^{i_{\mu}} \mathcal{C}^\mu_{\gd} \left( \pi_* (ev_{n+1}^* F_{\gd\mu} \otimes i_{\mu*}\ozzr)\right) .
   \end{align*}
   
     Here $f_i$ are polynomials with coefficients in $ev^*_{n+1}K^0(\ix )$, the bundles $E_\alpha, F_{\gd\mu}$ are elements of $K^0(\ix )$. To keep notation simple we write:
   \begin{align*}
   \Theta_{g,n,d}:= \AAA\cdot \BB \cdot \mathcal{C}_{g,n,d} .
    \end{align*}
  ``Twisted'' Gromov-Witten invariants are:
   \begin{align*}  
   \langle a_1\ops^{k_1},\ldots , a_n\ops^{k_n};\Theta \rangle_{g,n,d}: = \int_{[\ix_{g,n,d}]} \prod_{i=1}^n ev_i^*(a_i)\ops^{k_i}_i \cdot \Theta_{g,n,d} . 
    \end{align*}
Their generating series is the twisted potential $\mathcal{D}_{\mathcal{A,B,C}}$ :

    \begin{align*}
   & \mathcal{F}_{\mathcal{A,B,C}}^g (\mathbf{t}): = \sum_{d,n}\frac{Q^d}{n!}\langle \mathbf{t}(\ops),\ldots ,\mathbf{t}(\ops);\Theta  \rangle_{g,n,d}, \\
   &  \mathcal{D}_{\mathcal{A},\mathcal{B},\mathcal{C}}:= exp(\sum_g \hbar^{g-1} \mathcal{F}_{\mathcal{A,B,C}}^g ) .
    \end{align*}
  We view $ \mathcal{D}_{\mathcal{A},\mathcal{B},\mathcal{C}}$ as a formal function on $\mathcal{H}_+^{\mathcal{A,B,C}}$. 
  
   The symplectic vector space $(\mathcal{H}^{\mathcal{A,B,C}}, \Omega_{\mathcal{A,B,C}})$ is defined as $\mathcal{H}^{\mathcal{A,B,C}}=\mathcal{H}$, but with a different symplectic form :
   \begin{align*}
    \Omega_{\mathcal{A,B,C}}(\mathbf{f},\mathbf{g}):=\oint_{z=0} (\mathbf{f}(z),\mathbf{g}(-z))_{\mathcal{A}} dz 
    \end{align*} 
  where $($ , $)_{\mathcal{A}}$ is the {\em twisted} pairing given for $a,b\in H^*(I\ix)$ by:
     \begin{align*}
    (a,b)_{\mathcal{A}} := \langle a,b,1;\Theta\rangle_{0,3,0}.
     \end{align*}
       
  \begin{rem3}
  {\em We briefly discuss the case $(g,n,d)=(0,3,0)$. According to \cite{abgrvi} in this case the evaluation maps lift to $ev_i :\ix_{0,3,0}\to I\ix$. The spaces $\ix_{0,3,0,(\mu_1,\mu_2,0)}$ are empty unless $\mu_2 =\mu_1^I$, in which case they can be identified with $\ix_{\mu_1}$, with the evaluation maps being $ev_1 = id :\ix_{\mu_1}\to \ix_{\mu_1}$, $ev_2 = \iota : \ix_{\mu_1}\to \ix_{\mu^I_1}$ and $ev_3$ is the inclusion of $\ix_{\mu_1}$ in $\ix$.} 
  \end{rem3}
   \begin{altarem}
   {\em On $\ix_{0,3,0}$ there are no twistings of type $\mathcal{B}$ (the corresponding push-forwards are trivial for dimensional reasons) and of type $\mathcal{C}$ (there are no nodal curves). Hence the twisted pairing only depends on the $\mathcal{A}$ classes.}
   \end{altarem}  
  For a bundle $E$ on $\ix_\mu$ we denote by $E_{inv}$ the subbundle invariant under the action of the group element associated to $\ix_\mu$. According to the previous two remarks we can rewrite the pairing as:
   \begin{align*}
   (a,b)_{\mathcal{A}}:=\int_{I\ix}a\cdot\iota^*b \cdot \prod_\alpha \mathcal{A}_\ga\left((q^*E_\ga)_{inv}\right)  .
   \end{align*} 
 There is a rescaling map:
     \begin{align*}   
    (\mathcal{H}^{\mathcal{A},\mathcal{B},\mathcal{C}},\Omega_{\mathcal{A},\mathcal{B},\mathcal{C}})\to (\mathcal{H}, \Omega)\\
      a \mapsto a\sqrt{\prod_\alpha \mathcal{A}_\alpha((q^*E_\alpha)_{inv})}
     \end{align*}
   which identifies the symplectic spaces. We denote by $\mathcal{D}_{\mathcal{A,B}}, \mathcal{D}_{\mathcal{A}}$ etc. the potentials twisted only by classes of type occuring in the notation  and by
    \begin{align*} 
 [\ix_{g,n,d}]^{tw}:= [\ix_{g,n,d}]\cap \Theta_{g,n,d}  .
  \end{align*} 
  
     \section{Technical prerequisites}
   The computations in the proof of the theorems rely on pulling back the correlators on the universal orbicurve and on the locus of nodes.
Hence we need to know how the classes $\Theta_{g,n,d}$ behave under such pullbacks. The reader can skip this (unavoidably technical) section. To not make the statements and their proofs even more ugly we assume throughout this section that $i_\mu^{red}$ denotes the inclusion of a {\em single} nodal stratum  in the moduli space $\ix_{g,n+1,d}$. Otherwise equations $(\ref{29})$, $(\ref{restr1})$ and $(\ref{restr2})$ (and their proofs) need on the right hand side summation after all tuples $g_1+g_2=g$, $d_1+d_2=d$, $n_1+n_2=n$.  The result which we'll use in the proofs of the theorems is:     
        \begin{teta} \label{mainpllbc}
      The following equalities hold: 
 \begin{align} 
 1. \qquad &\pi^*[\ix_{g,n,d}]^{tw} = [\ix_{g,n+1,d}]^{tw} \cdot \prod_{\beta=1}^{i_B} \mathcal{B}_\beta \left( -\frac{f_\beta(L_{n+1}^{-1})-f_\beta(1)}{L_{n+1}-1}\right) +\nonumber\\
 &+\sum_{j=1}^n[\ix_{g,n+1,d}]^{tw}\cdot \left(\prod_{\delta=1}^{i_{\mu_j}}\mathcal{C}^{\mu_j}_\delta\left(- ev_{n+1}^*(F_{\delta\mu_j})\otimes \sigma_{j*}\mathcal{O}_{\mathcal{D}_j}\right) -1\right)+\nonumber\\
 &+ \sum_\mu[\ix_{g,n+1,d}]^{tw}\cdot\left( \prod_{\delta=1}^{i_\mu} \mathcal{C}^\mu_\delta\left(-ev_{n+1}^*(F_{\delta\mu})\otimes i_{\mu *}\ozzr \right)-1\right) .\label{301}\\
  2.\qquad &(\pi\circ i^{red}_\mu\circ p)^*[\ix_{g,n,d}]^{tw} = \nonumber\\
   & =  \frac{p_1^*([\ix_{g_1,n_1+1,d_1}]^{tw})\cdot p_2^*([\ix_{g_2,n_2+1,d_2}]^{tw})}{(ev_+^* \times ev_-^*)  \Delta_{\mu *}\prod_{\delta=1}^{i_{\mu}} \mathcal{C}^\mu_\delta\left( (q^*F_{\delta\mu})_\mu)\otimes (L_+L_--1)\right) } .\label{29}\\
   3.\qquad &(\pi\circ i^{irr}_{\mu}\circ p)^*[\ix_{g,n,d}]^{tw} = \nonumber\\
  & = \frac{[\ix_{g-1,n+2,d}]^{tw}}{ (ev_+^* \times ev_-^*)\Delta_{\mu *}\prod_{\delta=1}^{i_{\mu }} \mathcal{C}^\mu_\delta\left((q^*F_{\delta\mu})_\mu)\otimes (L_+L_--1)\right)}.\label{30}
 \end{align}
   \end{teta}
 Proof: all the equalities follow from the corresponding statements about the classes $\mathcal{A,B,C}$ separatedly, which we'll state and prove below. Formula $(\ref{301})$ follows from  $(\ref{aapi})$, $(\ref{bpp})$, $(\ref{nodalsum})$ combined with some more cancelation: namely the terms in $(\ref{nodalsum})$ supported on $\mathcal{D}_j$ and $\mathcal{Z}$ are killed by the correction factor in $(\ref{bpp})$ which is of the form $1+\psi_{n+1}\cdot ... $. The untwisted virtual fundamental classes satisfy $\pi^*[\ix_{g,n,d}] = [\ix_{g,n+1,d}]$. \\
   $(\ref{29})$ and $(\ref{30})$ follow from the corresponding Lemmata $\ref{tetaaaa}$, $\ref{tetabbb}$ and $\ref{tetaccc}$ for each of the classes $\AAA$, $\BB$ and $\CC$ combined with the splitting axiom in orbifold Gromov-Witten theory for the untwisted fundamental classes $[\ix_{g,n,d}]$, which we briefly review below. Let $\mathfrak{M}_{g,n}^{tw}$ be the stack of genus $g$ twisted curves with $n$ marked points. There is a natural map:
  \begin{align*}
   \mathrm{gl}:\mathfrak{D}^{tw}(g_1;n_1\vert g_2;n_2)\rightarrow \mathfrak{M}_{g,n}^{tw}
   \end{align*}
   induced by gluing two family of twisted curves into a reducible curve with a distinguished node. Here $\mathfrak{D}^{tw}(g_1;n_1\vert g_2,n_2)$ is defined as in  Section $5.1$ of \cite{abgrvi2}. This induces a cartesian diagram:

    \begin{displaymath}
     \begin{CD}
    \mathfrak{D}^{tw}_{g,n}(\ix )@>>>  \ix_{g,n,d} \\
      @VVV                 @VVV                          \\
 \mathfrak{D}^{tw}(g_1;n_1\vert g_2;n_2) @> \mathrm{gl} >> \mathfrak{M}_{g,n}^{tw}.
             \end{CD}
 \end{displaymath}
 There is a natural map:
  \begin{align*}
 \mathfrak{g}:\bigcup_{d_1+d_2=d} \ix_{g_1,n_1+1,d_1}\times_{\overline{I\ix}} \ix_{g_2,n_2+1,d_2}\rightarrow \mathfrak{D}^{tw}_{g,n}(\ix).
  \end{align*}
 Then the diagram:

    \begin{displaymath}
     \begin{CD}
     \ix_{g_1,n_1+1,d_1}\times_{\overline{I\ix}} \ix_{g_2,n_2+1,d_2} \subset \mathcal{Z}@>>>  \overline{I\ix}\\
      @VVV                 @V  \Delta VV                          \\
 \ix_{g_1,n_1+1,d_1}\times \ix_{g_2,n_2+1,d_2} @>ev_+ \times \breve{ev}_- >> \overline{I\ix}\times \overline{I\ix}
             \end{CD}
 \end{displaymath}
 gives:
  \begin{align}
  \sum_{d_1+d_2 =d}\Delta^!([\ix_{g_1,n_1+1,d_1}]\times [\ix_{g_2,n_2+1,d_2}]) =\mathfrak{g}^*(\mathrm{gl}^!([\ix_{g,n,d}])).\label{46}
   \end{align}
  For details and proofs of the statements we refer the reader to the paper \cite{abgrvi2} (Prop. $5.3.1.$). The only modification we have made is - we consider the class of the diagonal with respect to the \emph{twisted} pairing on $\ix_{0,3,0,(\mu_1,\mu_2,0)}$. This cancels the factor $ev_\Delta^*(\mathcal{A}_{0,3,0})$ in $(\ref{29})$ and $(\ref{30})$.
 
    Informally relation $(\ref{46})$ says that the restriction of the virtual fundamental class of $\ix_{g,n,d}$ to $\mathcal{Z}$ coincides with the push forward  of the product of virtual fundamental classes under the gluing morphisms. Hence integration on $\mathcal{Z}$ factors "nicely" as products of integrals on the two separate moduli spaces. \vspace{10pt}
   
   The rest of the section is devoted to proving pullback results about each type of twisting class separately.  
     \begin{pushpull}\label{sss}
 Consider the following diagram:
       \begin{displaymath}
     \begin{CD}
    \mathcal{X}_{g,n+\circ+\bullet,d,(\mu_1,\ldots ,\mu_n,0,0)} @>\pi_1 >>  \mathcal{X}_{g,n+\bullet , d,(\mu_1,\ldots ,\mu_n,0)}\\
      @V \pi_2 VV                 @V  \pi_2 VV                          \\
      \mathcal{X}_{g,n+\circ,d,(\mu_1,\ldots ,\mu_n,0)} @>\pi_1 >>  \ix_{g,n,d,(\mu_1,\ldots , \mu_n)}
             \end{CD}
 \end{displaymath}
where $\pi_1$ forgets the $(n+1)$-st marked point (which we denoted $\circ$) and $\pi_2$ forgets the $(n+2)$-nd marked point (denoted $\bullet$)
and let $\alpha\in K^0(\mathcal{X}_{g,n+\circ,d,(\mu_1,\ldots ,\mu_n,0)})$. Then $\pi_2^*\pi_{1*}\alpha = \pi_{1*}\pi_2^*\alpha$.
    \end{pushpull}

Proof: for simplicity of notation we suppress the labeling $(\mu_1, \ldots , \mu_n)$ in the proof. Consider the fiber product:
   \begin{align*}
     \mathcal{F} :=  \mathcal{X}_{g,n+\circ,d}\times_{ \ix_{g,n,d}} \mathcal{X}_{g,n+\bullet , d }
    \end{align*}
  and denote by $p_1, p_2$ the projections from $\mathcal{F}$ to the factors and by $\varphi: \mathcal{X}_{g,n+\circ+\bullet,d}\rightarrow \mathcal{F}$ the morphism induced by $\pi_1 , \pi_2$.  $\varphi$ is a birational map: it has positive dimensional fibers along the locus where the two extra marked points hit another marked point or a node. This locus has codimension $2$ - this in particular shows that $\mathcal{F}$ is normal. We'll prove that  
        \begin{align*}
    \varphi_*(\mathcal{O}_{\mathcal{X}_{g,n+\circ+\bullet,d}})= \mathcal{O}_{\mathcal{F}}. 
    \end{align*}
    By definition of $K$-theoretic push-forward
     \begin{align*}
     \varphi_*\mathcal{O}_{\mathcal{X}_{g,n+\circ+\bullet,d}} = R^0 \varphi_* \mathcal{O}_{\mathcal{X}_{g,n+\circ+\bullet ,d}} -  R^1 \varphi_*\mathcal{O}_{\mathcal{X}_{g,n+\circ+\bullet,d}}.
      \end{align*}
  It is easy to see that $R^0\varphi_*(\mathcal{O}_{\mathcal{X}_{g,n+\circ+\bullet,d}})= \mathcal{O}_{\mathcal{F}}$ as quasicoherent sheaves (this is true for every proper birational map with normal target). We only have to prove that $R^1 = 0$, which we do by looking at the stalks:
    \begin{align*}
     (R^1 \varphi_* \mathcal{O}_{\mathcal{X}_{g,n+ \circ +\bullet , d}})_x  = H^1( \varphi^{-1}(x), \mathcal{O}_{\mathcal{X}_{g,n+\circ+ \bullet ,d}\vert \varphi^{-1}(x)}).
      \end{align*}
      If the fiber over $x$ is a point, there's nothing to prove. If $x$ is in the blowup locus the fiber is a (possibly weighted) $\mathbb{P}^1$. A calculation in \cite{abgrvi2} shows that :
   \begin{align*}
   \chi (\mathcal{C},\mathcal{O}_{\mathcal{C}} ) = 1-g ,
   \end{align*} 
    where $g$ is the arithmetic genus of the coarse curve $C$. this shows that $H^1(\varphi^{-1}(x),\mathcal{O} )=0$ .
  We have $p_{1*}p_2^*\alpha = \pi_2^*\pi_{1*}\alpha $ because the diagram:
      \begin{displaymath}
     \begin{CD}
    \mathcal{F} @>p_1 >>  \mathcal{X}_{g,n+\bullet , d,(\mu_1,\ldots ,\mu_n,0)}\\
      @V p_2 VV                 @V  \pi_2 VV                          \\
      \mathcal{X}_{g,n+\circ,d,(\mu_1,\ldots ,\mu_n,0)} @>\pi_1 >>  \ix_{g,n,d,(\mu_1,\ldots , \mu_n)}
             \end{CD}
 \end{displaymath}
is a fiber square. Therefore:
 \begin{align*}
  \pi_{1*}\pi_2^*\alpha = p_{1*}\varphi_* \left(\varphi^*p_2^*\alpha\right)  = p_{1*}p_2^*\alpha \varphi_*(\mathcal{O})  = p_{1*}p_2^*\alpha = \pi_2^*\pi_{1*}\alpha .
  \end{align*}

  We  need to know how the classes $\AAA , \BB , \CC $ behave under pullback by the morphisms $\pi$ and $\pi\circ i\circ p$.
 
    \begin{tetaaa}  \label{tetaaaa}
 The following identities hold:
   \begin{align}
  a.\quad &\pi^* \AAA = \mathcal{A}_{g,n+1,d}.   \label{aapi}\\
  b.\quad &(\pi\circ i^{red}\circ p )^* \AAA = \frac{ p_1^* \mathcal{A}_{g_1,n_1+1,d_1} \cdot p_2^* \mathcal{A}_{g_2,n_2 +1,d_2}}{ev_\Delta^*\mathcal{A}_{0,3,0}}.\label{restr1}\\
  c.\quad & (\pi\circ i^{irr}\circ p)^* \AAA = \frac{\mathcal{A}_{g-1,n+2,d}}{ev_\Delta^*\mathcal{A}_{0,3,0}}.
   \end{align}
 \end{tetaaa}
  Denote by $E_{g,n,d}:= \pi_*(ev_{n+1}^* E)$. Then it is shown in \cite{tseng} that:
  \begin{align*}
  a.\quad & \pi^* E_{g,n,d} =  E_{g,n+1,d},\\
  b.\quad &(\pi\circ i^{red}\circ p)^* E_{g,n,d} =p_1^*(E_{g_1,n_1+1,d_1}) + p_2^*(E_{g_2,n_2+1,d_2})-ev^*_\Delta (q^* E_{inv}),\\
  c.\quad & (\pi\circ i^{irr}\circ p)^* E_{g,n,d} = E_{g-1,n+2,d} - ev_\Delta^*(q^*E_{inv}).
   \end{align*}
  The identities then follow by multiplicativity of the classes $\mathcal{A}_\alpha$. We regard the class $\mathcal{A}_{0,3,0}$
 as an element of $H^*(I\ix ,\mathbb{Q})$. We can then pull it back by the diagonal evaluation morphism $ev_\Delta$ at the node.    
 
   \begin{tetabb} \label{tetabbb} The following hold:
    \begin{align}
   a.\quad  & \pi^* \BB = \mathcal{B}_{g,n+1,d}\cdot  \prod_{\beta=1}^{i_B} \mathcal{B}_\beta \left( -\frac{f_\beta(L_{n+1}^{-1})-f_\beta(1)}{L_{n+1}-1}\right). \label{bpp}\\
   b.\quad  & (\pi\circ i^{red}\circ p)^* \BB =  p_1^* \mathcal{B}_{g_1,n_1+1,d_1} \cdot p_2^* \mathcal{B}_{g_2,n_2+1,d_2}. \label{restr2}\\
   c.\quad  & (\pi\circ i^{irr}\circ p)^* \BB =  \mathcal{B}_{g-1,n+2,d}. \label{10}
     \end{align}
    \end{tetabb}
  Proof: The first identity is a consequence of Lemma $\ref{sss}$ . More precisely we apply the lemma to the class $\alpha = ev_{n+1}^*(E)(L_{n+1}-1)^{k+1} $. This gives:
  \begin{align*}
 & \pi_2^*\pi_{1*} \left[ ev_{n+1}^*(E)(L_{n+1}-1)^{k+1}\right] = \pi_{1*}\pi_2^* \left[ ev_{n+1}^*(E)(L_{n+1}-1)^{k+1}\right] =\\
 & = \pi_{1*} \left[ ev_{n+1}^*(E) (L_{n+1}-1)^{k+1}- (\sigma_{\bullet})_*\left(ev_{n+1}^*(E)(L_{n+1}-1)^{k}\right)\right] =\\
 &  = \pi_{1*} \left( ev_{n+1}^*(E) (L_{n+1}-1)^{k+1}\right) - ev_{n+1}^*(E)(L_{n+1}-1)^{k}.
   \end{align*}
  The last equality follows because $\pi_1\circ \sigma_\bullet = Id$ and the second equality uses the comparison identity for cotangent line bundles $L_i$:
    \begin{align*}
   \pi^*( (L_{i}-1)^{k+1}) = (L_{i}-1)^{k+1} - \sigma_{i*}\left[(L_{i}-1)^{k}\right].
     \end{align*}
But both morphisms $\pi_1, \pi_2$ can be identified with the universal orbicurve $\pi$. Hence we deduce:
 \begin{align}
\pi^*\pi_{*} \left( ev_{n+1}^*(E)(L_{n+1}-1)^{k+1}\right)& = \pi_{*} \left( ev_{n+2}^*(E) (L_{n+2}-1)^{k+1}\right) - \nonumber\\
   &-ev_{n+1}^*(E)(L_{n+1}-1)^{k} ,\label{61}
 \end{align}
 or more generally if we expand 
  \begin{align*}
   f_\beta(L_{n+1}^{-1})-f_\beta(1)=\sum_{k\geq 0} a_k (L_{n+1}-1)^{k+1},
   \end{align*}
    then:
 \begin{align}
\pi^*\pi_{*} (f_\beta(L_{n+1}^{-1})-f_i(1)) = \pi_*(f_\beta(L_{n+2}^{-1})-f_i(1)) - \frac{f_\beta(L_{n+1}^{-1})-f_\beta(1)}{L_{n+1}-1}. \label{62}
  \end{align}
  Then $(\ref{bpp})$ follows because $\mathcal{B}_\beta$ are multiplicative classes:
   \begin{align*}
 \pi^*\mathcal{B}_\beta &\left(\pi_*(f_\beta(L_{n+1}^{-1})-f_\beta(1))\right) =  \mathcal{B}_\beta\left(\pi^*\pi_*(f_\beta(L_{n+1}^{-1})-f_\beta(1))\right)= \\
                 & = \mathcal{B}_\beta\left(\pi_*(f_\beta(L_{n+2}^{-1})-f_\beta(1))- \frac{f_\beta(L_{n+1}^{-1})-f_\beta(1)}{L_{n+1}-1}\right) = \\
              & = \mathcal{B}_\beta\left( \pi_*(f_\beta(L_{n+2}^{-1})-f_\beta(1))\right)\cdot \mathcal{B}_\beta\left(- \frac{f_\beta(L_{n+1}^{-1})-f_\beta(1)}{L_{n+1}-1}\right).
    \end{align*}
  
  \begin{exbet} 
  \emph{ In the case $f_\beta = ev_{n+1}^* (E_\beta) \otimes L_{n+1}^{-1}$ (which is the only one we'll need) we have:}
   \begin{align*}
  \frac{f_\beta(L_{n+1}^{-1})-f_\beta(1)}{L_{n+1}-1} = -E_\beta L_{n+1}^{-1}
   \end{align*}
  \emph{and relation $(\ref{bpp})$ reads:}
\begin{align}
\pi^* \BB = \mathcal{B}_{g,n+1,d}\cdot  \prod_{\beta=1}^{i_B}\mathcal{B}_\beta( E_\beta\otimes L_{n+1}^{-1}).
\end{align}
  \end{exbet}
 Relation  $(\ref{restr2})$ follows from the identity:
   \begin{align*}
   (\pi \circ i_{red})^*[ \pi_*(f(L_{n+1}^{-1})& - f(1))] = \\
      =  & p_1^*[\pi_*(f(L_{n_1+2}^{-1})- f(1))] + p_2^*[ \pi_*(f(L_{n_2+2}^{-1})- f(1))],
   \end{align*}
  which we prove below. By linearity is enough to prove the result for $f=(L_{n+1}-1)^{k+1}$ for $k\geq 0$. Assume for now that $k\geq 1$. Relation $(\ref{61})$ gives:
  \begin{align}
  \pi^*\pi_{*} (L_{n+1}-1)^{k+1} = \pi_{*}(L_{n+2}-1)^{k+1}- (L_{n+1}-1)^k . \label{12332}
  \end{align}
  When we apply $p^* i_{red}^*$ to this relation the second summand in the RHS of $(\ref{12332})$ vanishes because $L_{n+1}$ is trivial on $\widetilde{\mathcal{Z}}$. Therefore
 \begin{align*}
p^* i_{red}^*\pi^*\pi_{*} (L_{n+1}-1)^{k+1} = (i_{red}\circ p)^* \pi_{*} (L_{n+2}-1)^{k+1}.
 \end{align*}
 Let $\ix_{g_1,n_1+1,d_1}\times_{\overline{I\ix}}\ix_{0,3,0}\times_{\overline{I\ix}}\ix_{g_2,n_2+1,d_2}$  be a stratum of $\mathcal{Z}$. If we denote by $\pi: \mathcal{U}_{g,n,d}'\rightarrow \mathcal{U}_{g,n,d}$ the universal curve then we have a fiber diagram: 
   \begin{displaymath}
     \begin{CD}
    \mathcal{Z}_1\cup  \mathcal{Z}_2 \cup \mathcal{Z}_3 @>i >>  \mathcal{U}_{g,n,d}'\\
      @V \pi VV                 @V  \pi VV                          \\
    \ix_{g_1,n_1+1,d_1}\times_{\overline{I\ix}}\ix_{0,3,0}\times_{\overline{I\ix}}\ix_{g_2,n_2+1,d_2}  @> i >>  \mathcal{U}_{g,n,d}  .
             \end{CD}
 \end{displaymath}
Here $\mathcal{Z}_1$ and $\mathcal{Z}_3$ are the universal curves over the factors $\ix_{g_1,n_1+1,d_1}$ and $\ix_{g_2,n_2+1,d_2}$. So using
  \begin{align}
  i_{red}^* \pi_{*} (L_{n+2}-1)^{k+1} =\pi_* i_{red}^* (L_{n+2}-1)^{k+1} ,
   \end{align}
   we see that the contribution of the strata $\mathcal{Z}_1$ and $\mathcal{Z}_3$ above is:
   \begin{align}
   p_1^*[\pi_*(f(L_{n_1+2}^{-1})- f(1))] + p_2^*[ \pi_*(f(L_{n_2+2}^{-1})- f(1))].
   \end{align}
  So if we show that the contribution from $\mathcal{Z}_2$ is $0$ we are done.$\mathcal{Z}_2$ is the universal curve over the factor $\ix_{0,3,0}$, hence it is a fiber product $\ix_{g_1,n_1+1,d_1}\times_{\overline{I\ix}}\ix_{0,4,0}\times_{\overline{I\ix}}\ix_{g_2,n_2+1,d_2}$ . The fibers of the map $\mathcal{Z}_2\rightarrow \mathcal{Z}$ are (weighted) $\mathbb{P}^1$. However the class $L_{n+2}$ (consider it as the cotangent line $L_1\in K^0(\overline{M}_{0,4})$) is a cotangent line at a point with trivial orbifold structure, so we can use Y. P. Lee's formula in \cite{ypl} which in this particular case reads:
   \begin{align}
     \chi(\overline{M}_{0,4}, L_1^k ) = k+1 .
    \end{align} 
   Hence the Euler characteristics of $(L_{n+2} - 1)^{k +1}$ is: 
    \begin{align*}
    &\chi\left(\overline{M}_{0,4}, (L_{1} - 1)^{k +1}\right) = \sum^{k+1}_{i=0} (i+1)(-1)^{k+1-i} \binom{k+1}{i} = \nonumber\\
     &  =  \sum^{k+1}_{i=0}(-1)^{k+1-i} \binom{k+1}{i} + (k+1)\sum^{k+1}_{i=1}(-1)^{k+1-i} \binom{k}{i-1} = 0+0 = 0 .
    \end{align*}
  This almost proves the statement. We are left with the case $k=0$, which is slightly different: the sum above equals $1$, but this is cancelled by the $-1$ in the second term of $(\ref{12332})$. Relation $(\ref{restr2})$ follows then from the multiplicativity of the classes $\mathcal{B}_\beta$.
    A similar computation shows relation $(\ref{10})$.
   
  \begin{chernch} \label{nodp5}
 Let $F\in K^0(\ix)$. Then:
  \begin{align}
 a.\quad & \pi^*\pi_{*} i_{\mu*} (ev^*_{n+1}(F) \otimes \ozzr)  = \pi_{*} i_{\mu*} (ev^*_{n+1}(F) \otimes \ozzr )  - \nonumber\\
   & \sum_{j,\mu_j=\mu}  ev^*_{n+1}(F) \otimes \sigma_{j*}\mathcal{O}_{\mathcal{D}_j} -    i_{\mu*} (ev^*_{n+1}(F)\otimes \ozzr ) . \label{nodp3}\\
 b.\quad & (\pi\circ i\circ p)^*(\pi_{*} i_{\mu*} (ev^*_{n+1}(F) \otimes \ozzr ) =  p_1^* (\pi_{*} i_{\mu*} (ev^*_{n+1}(F) \otimes \ozzr ) ) + \nonumber\\
  &   p_2^*(\pi_{*} i_{\mu*} (ev^*_{n+1}(F) \otimes \ozzr ))+  \left( ev_{n+1}^* F \otimes(1-L_+L_- )\right)  .\label{nodp4}
  \end{align}
    \end{chernch}
   \begin{rem11}
  \emph{Before delving in the technicalities of the proof, we try a heuristic explanation of why the rather ugly formulae ``should'' be true:}
      \begin{itemize}
     \item {\em Assume for now that $F$ is the trivial bundle $\mathbb{C}$. The nodal locus $\mathcal{Z}$ ``separates nodes'' in the following sense: above a point of $\ix_{g,n,d}$ representing a nodal curve with $k$ nodes lie exactly $k$ points of $\mathcal{Z}$. This is very similar with the way the normalization of a nodal curve $ \widetilde{C}\to C$ separates the nodes. But the structure sheaves of $\widetilde{C}$ and $C$ differ (in K-theory) by skyscraper sheaves at the preimages of nodes. That's pretty much what the first formula expresses: the pull-back of the structure sheaf of the codimension one stratum of nodal curves in $\ix_{g,n,d}$ equals the structure sheaf of the nodal locus in the universal family, minus a copy of the structure sheaf of $\mathcal{Z}$ (which has codimension two in the universal family) itself. The terms supported on the divisors $\mathcal{D}_j$ are substracted because they are nodes in the universal family, but they lie over the whole space $\ix_{g,n,d}$. We'll see that the presence of the class $ev^*_{n+1}(F)$ doesn't complicate things too much.}
     \item {\em For the second formula, think of $\pi_* i_{\mu *}\alpha$ as a class supported on  a codimension one subvariety. We pull it back along the map ($\pi i$)
, which is like restricting to another codimension one subvariety. If these subvariesties intersect along a codimension two cycle (represented by curves with two nodes), then they contribute $p^*_i(\pi_* i_{\mu*}\alpha )$ to $(\ref{nodp4})$. If they are the same subvariety, then $\alpha$ gets multiplied with the Euler class of the normal bundle of it in the ambient space, which is $1-L_+L_-$.} 
      \end{itemize}   
   \end{rem11} 
 
     Proof of Lemma $\ref{nodp5}$: denote by $\mathcal{Z}_\bullet$, $\mathcal{Z}_{\circ}$, respectively $\mathcal{Z}_{\bullet\circ}$ the nodal loci living inside the corresponding moduli spaces (and by $\mathcal{Z}_{\circ,\mu}$ etc. the ones with nodes of specific orbifold type) in the following diagram (we write $\overline{\mu}$ for the sequence $(\mu_1,\ldots ,\mu_n)$): 
     \begin{displaymath}
     \begin{CD}
     \pi_2^{-1}(\mathcal{Z}_{\circ,\mu}) @>i_\mu >>\cup_{\overline{\mu}} \mathcal{X}_{g,n+\circ+\bullet,d,(\overline{\mu},0,0)} @>\pi_1 >> \cup_{\overline{\mu}} \mathcal{X}_{g,n+\bullet, d,(\overline{\mu},0)}\\
      @V \pi_2 VV                                    @V  \pi_2 VV              @V \pi_2 VV              \\
   \mathcal{Z}_{\circ,\mu} @>i_\mu >>   \cup_{\overline{\mu}} \mathcal{X}_{g,n+\circ,d,(\overline{\mu},0)} @>\pi_1 >>  \ix_{g,n,d}  .
             \end{CD}
      \end{displaymath}

 Remember that $\mathcal{Z}_{\circ,\mu}$ is defined as the total range of the gluing map (for simplicity we omit in the notation the stratum parametrizing self-intersecting curves; the proof carries through word by word):
     \begin{align*}
  \ix_{g_1,n_1+1,d_1}\times_{\overline{\ix}_\mu\times \overline{\ix}_{\mu^I}} \ix_{0,3,0}\times_{\overline{\ix}_\mu\times \overline{\ix}_{\mu^I}}\ix_{g_2,n_2+1,d_2}\rightarrow\mathcal{Z}_\circ\hookrightarrow  \ix_{g,n+\circ ,d}.
     \end{align*}
  We will compute  $\pi_2^*(\pi_{1*} i_{\mu*} (ev^*_{\circ}(F) \otimes \mathcal{O}_{\mathcal{Z}_{\circ,\mu}}))$.      
     
The square on the left is a fiber diagram, hence $i_*\pi_{2}^* =\pi_{2}^*i_*$. For the one on the right we have proved that  $\pi_2^*\pi_{1*}= \pi_{1*}\pi_{2}^*$. Therefore:
   \begin{align}
    \pi_2^*(\pi_{1*} i_{\mu*} (ev^*_{n+1}(F) \otimes \mathcal{O}_{\mathcal{Z}_{\circ,\mu}}))= \pi_{1*} i_{\mu*}\pi_2^* (ev^*_{\circ}(F) \otimes \mathcal{O}_{\mathcal{Z}_{\circ,\mu}})) . \label{7676}
   \end{align} 
  But:
   \begin{align*}
   \pi_2^* (ev^*_{\circ}F \otimes \mathcal{O}_{\mathcal{Z}_{\circ,\mu}})=ev^*_{\circ}F\otimes\mathcal{O}_{\pi^{-1}_2(\mathcal{Z}_{\circ,\mu})} .
   \end{align*}
 
 The space $\pi^{-1}_2(\mathcal{Z}_{\circ,\mu}): = \mathcal{Z}_{\circ,1}\cup \mathcal{Z}_{\circ,2}\cup \mathcal{Z}_{\circ,3}$ is a singular space, where each codimension two stratum is the universal curve over one factor of $\mathcal{Z}_{\circ,\mu}$ and they intersect along two codimension three strata, call them $\mathcal{Z}_{12}$ and $\mathcal{Z}_{23}$:
   \begin{align*}
\mathcal{Z}_{12}=\ix_{g_1,n_1+1,d_1} \times_{\overline{I\ix}} \ix_{0,3,0}\times_{\overline{I\ix}} \ix_{0,3,0}\times_{\overline{I\ix}} \ix_{g_2,n_2+1,d_2} 
   \end{align*}
  where the two rational components carry the points $\bullet$, $\circ$ and two nodes. Figure $1$ schematically represents each of these five strata.
    \begin{figure}
\begin{center}
\scalebox{1.0}{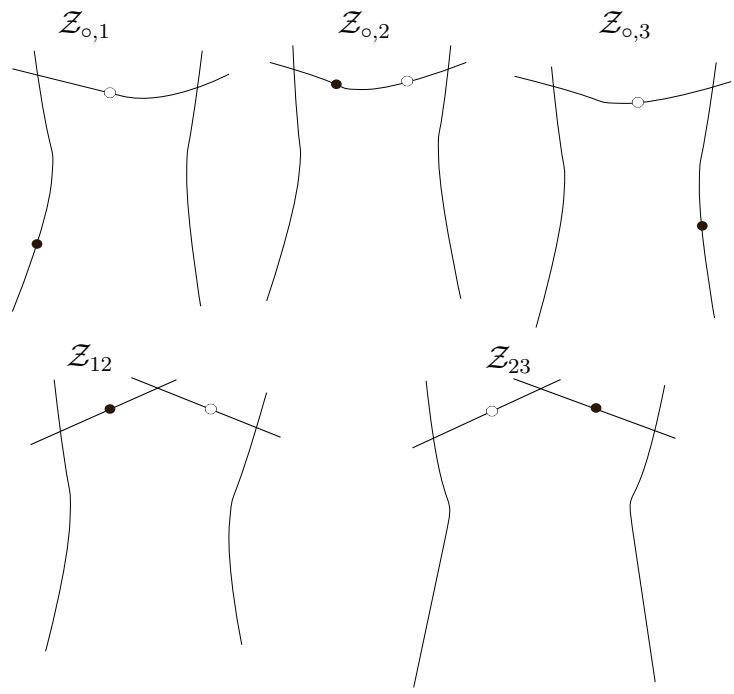}
\caption{Strata of $\pi^{-1}_2(\mathcal{Z}_{\circ,\mu})$.}
\label{fig:latexexport}
\end{center}
\end{figure} 
  We can write the structure sheaf of $\pi^{-1}_2(\mathcal{Z}_{\circ,\mu})$ as:
  \begin{align*}   
  \mathcal{O}_{\pi^{-1}_2(\mathcal{Z}_{\circ,\mu})} =\mathcal{O}_{\mathcal{Z}_{\circ,1}} + \mathcal{O}_{\mathcal{Z}_{\circ,3}} + \mathcal{O}_{\mathcal{Z}_{\circ,2}} - \mathcal{O}_{\mathcal{Z}_{12}} - \mathcal{O}_{\mathcal{Z}_{23}} .
  \end{align*}  
   
  We tensor this with the class $ev_{\circ}^* F$, keeping in mind that on the strata $\mathcal{Z}_{\circ,2},\mathcal{Z}_{12},\mathcal{Z}_{23}$ $ev_\circ = ev_\bullet$:
  \begin{align}   
 ev_\circ^*F \otimes\mathcal{O}_{\pi^{-1}_2(\mathcal{Z}_{\circ,\mu})} =  ev_\circ^*F\otimes\left[\mathcal{O}_{\mathcal{Z}_{\circ,1}} + \mathcal{O}_{\mathcal{Z}_{\circ,3}}\right] +    ev_\bullet^*F\otimes\left[\mathcal{O}_{\mathcal{Z}_{\circ,2}} - \mathcal{O}_{\mathcal{Z}_{12}} -  \mathcal{O}_{\mathcal{Z}_{23}}\right] .\label{7778}
  \end{align}
   We plug $(\ref{7778})$ in $(\ref{7676})$ and we get:  
 \begin{align}
    & \pi_2^*(\pi_{1*} i_{\mu*} (ev^*_{\circ}(F) \otimes \mathcal{O}_{\mathcal{Z}_{\circ,\mu}}))= \nonumber\\
    & = \pi_{1*} i_{\mu*}\left[ ev_\circ^*F\left(\mathcal{O}_{\mathcal{Z}_{\circ,1}} +\mathcal{O}_{\mathcal{Z}_{\circ,3}}\right) + 
  ev_\bullet^*F\left(\mathcal{O}_{\mathcal{Z}_{\circ,2}} -\mathcal{O}_{\mathcal{Z}_{12}} -  \mathcal{O}_{\mathcal{Z}_{23}}\right) \right].\label{7780}
    \end{align} 
    
  We now notice that the union of $\mathcal{Z}_{\circ,1}$ and $\mathcal{Z}_{\circ,3}$ is almost $\mathcal{Z}_{\bullet\circ,\mu}$ but not quite. 
 There are strata:
   \begin{align*}
 \ix_{g,n,d}\times_{\overline{\ix}_\mu} \ix_{0,3,0} \times_{\overline{\ix}_\mu} \ix_{0,3,0} 
  \end{align*}
  which are in $\mathcal{Z}_{\bullet\circ,\mu}$, but they are missing from $\mathcal{Z}_{\circ,1} \cup \mathcal{Z}_{\circ,3}$ because the map $\pi_2\circ i_\mu$ contracts one rational tail. These are mapped by $\pi_1\circ i_\mu$ isomorphically to divisors $\mathcal{D}_j\in \ix_{g,n+\bullet,d}$. There is one such stratum for each $j$ such that $\mu_j = \mu$. Hence we can write:
  \begin{align}
\pi_{1*}i_{\mu*} \left[ev_\circ^*F\mathcal{O}_{\mathcal{Z}_{\circ,1}} +  ev_\circ^*F\mathcal{O}_{\mathcal{Z}_{\circ,3}}\right] = 
   \pi_{1*} i_{\mu*} (ev^*_{\circ}(F) \otimes \ozzr )  - \nonumber\\
   - \sum_{j,\mu_j=\mu}  ev^*_{\circ}(F) \otimes \sigma_{j*}\mathcal{O}_{\mathcal{D}_j}. \label{7777}
  \end{align}  
 The codimension three strata $\mathcal{Z}_{12}$ and $\mathcal{Z}_{23}$ are mapped by $\pi_{1}i_{\mu}$ isomorphically  to $\mathcal{Z}_{\bullet,\mu}$. As for $\mathcal{Z}_{\circ,2}$, this is a $\mathbb{P}^1$ fibration over $\mathcal{Z}_{\bullet,\mu}$. When we push forward, we integrate the structure sheaf of (weighted) $\mathbb{P}^1$. This equals $1$, as already explained. At the end of the day we see that the last three terms in $(\ref{7780})$ contribute:
\begin{align}
 \pi_{1*}i_{\mu*}\left[ev_\bullet^*F\left(\mathcal{O}_{\mathcal{Z}_{\circ,2}} - \mathcal{O}_{\mathcal{Z}_{12}} -  \mathcal{O}_{\mathcal{Z}_{23}}\right)\right] = -ev^*_\bullet F\otimes i_{\mu_*}\mathcal{O}_{\mathcal{Z}_{\bullet,\mu}}.  \label{7799}
  \end{align}
   Adding up $(\ref{7777})$ with $(\ref{7799})$ and identifying $\pi_1=\pi_2 =\pi$ and $ev_\circ =ev_{n+1}$ proves the first  equality in the lemma.
   
  For the second equality, we first prove: 
      \begin{nodalclass4} \label{lemma11}
  Let $j:\mathcal{Z}\hookrightarrow \mathcal{U}_{g,n,d}$ be the codimension two nodal locus. Then:
  \begin{align}
 & j^* \pi_*i_{\mu*}\left( ev_{n+1}^* F \otimes \ozzr \right) = p_1^* \pi_*i_{\mu*}\left( ev_{n+1}^* F \otimes \ozzr\right)+  \nonumber\\
   & + p^*_2 \pi_*i_{\mu*}\left( ev_{n+1}^* F \otimes \ozzr \right) + (2-L_+ - L_- ) ev^*_{n+1}(F).  \label{500}
  \end{align} 
 \end{nodalclass4}
  Proof of the Lemma $\ref{lemma11}$: let $\mathcal{U}_{g,n,d}'$ be the universal curve over $\mathcal{U}_{g,n,d}$. The universal curve over $\mathcal{Z}$ is a union of three types of strata, depending on which component the extra marked point on $\mathcal{U}_{g,n,d}'$ - which we denote $\bullet$ - lies on :
     \begin{align*}
    \mathcal{Z}_1 = \ix_{g_1,n_1+1+\bullet,d_1}\times_{\overline{I\ix}}\ix_{0,3,0}\times_{\overline{I\ix}} \ix_{g_2,n_2+1,d_2},\\  
     \mathcal{Z}_2 = \ix_{g_1,n_1+1,d_1}\times_{\overline{I\ix}}\ix_{0,3+\bullet,0}\times_{\overline{I\ix}} \ix_{g_2,n_2+1,d_2},\\
      \mathcal{Z}_3 = \ix_{g_1,n_1+1,d_1}\times_{\overline{I\ix}}\ix_{0,3,0}\times_{\overline{I\ix}} \ix_{g_2,n_2+1+\bullet,d_2}.
     \end{align*}
      The diagram below is a fiber square:
    \begin{displaymath}
     \begin{CD}
    \mathcal{Z}_1\cup  \mathcal{Z}_2 \cup \mathcal{Z}_3 @>j >>  \mathcal{U}_{g,n,d}'\\
      @V \pi VV                 @V  \pi VV                          \\
  \mathcal{Z} @> j >>  \mathcal{U}_{g,n,d}.
             \end{CD}
 \end{displaymath}
    Hence : $j^*\pi_* i_{\mu*}\ga =\pi_* j^*i_{\mu*} \ga$. To compute $j^*i_{\mu*}\ga$ we form the following fiber diagram: 
   \begin{displaymath} 
     \begin{CD}
      \overline{\mathcal{Z}} @> j >>  \mathcal{Z}_{\bullet,\mu} \\
      @V \pi VV                 @V  \pi VV                          \\
    \mathcal{Z}_1\cup  \mathcal{Z}_2 \cup \mathcal{Z}_3 @>j >>  \mathcal{U}_{g,n,d}' .
      \end{CD}
 \end{displaymath}
    The space $\overline{\mathcal{Z}}$ is simply the intersection of $\mathcal{Z}_1\cup \mathcal{Z}_2\cup\mathcal{Z}_3$ with $\mathcal{Z}_{\bullet,\mu}$. Where the intersection is transversal one can simply write $j^*i_{\mu*}\ga =i_{\mu*}j^*\ga$. On components where the intersection is not transversal, there is some excess bundle $N$ and $j^*i_{\mu*}\ga =i_{\mu*} e(N)j^*\ga$. The strata $\mathcal{Z}_1$ and $\mathcal{Z}_3$ intersect the nodal locus $\mathcal{Z}_{\bullet,\mu}$ in $\mathcal{U}_{g,n,d}'$ transversely along  codimension four strata which can be seen as the nodal locus in $\ix_{g_1,n_1+1+\bullet,d_1}$ and $\ix_{g_2,n_2+1+\bullet,d_2}$ respectively. Hence the contribution to $(\ref{500})$ is:
     \begin{align*}
  p_1^* \pi_*i_{\mu*}\left( ev_{n+1}^* F \otimes \ozzr \right) 
    + p^*_2\pi_*i_{\mu*}\left( ev_{n+1}^* F \otimes \ozzr\right).
      \end{align*}
   On the other hand $\mathcal{Z}_2$ intersects $\mathcal{Z}_{\bullet,\mu}$ along two codimension three strata of the form:
   \begin{align*}
  \mathcal{Z}_1 = \ix_{g_1,n_1+1,d_1}\times_{\overline{I\ix}}\ix_{0,3,0}\times_{\overline{I\ix}} \ix_{0,3,0}\times_{\overline{I\ix}}\ix_{g_2,n_2+1,d_2}.\\  
   \end{align*}   
   Each gives a one dimensional excess normal bundle with Euler classes $1-L_+$ and $1-L_-$ respectively. They project isomorphically to $\mathcal{Z}$ downstairs. Hence they contribute:
   \begin{align*}
   (2-L_+ -L_-) ev^*_{n+1}(F).
   \end{align*}
    Adding up, we get $(\ref{500})$.
  
   We now prove formula $\ref{nodp4}$ in Lemma $\ref{nodp5}$. It falls out easily by combining $(\ref{nodp3})$ with  Lemma $\ref{lemma11}$. More precisely we take $i^*$ of formula $(\ref{nodp3})$: the first term is computed in Lemma $\ref{lemma11}$, the part supported on $\mathcal{D}_j$ vanishes and:
    \begin{align}
    i_\mu^*i_{\mu*}\ozzr = e(N) = (1-L_-)(1-L_+)
     \end{align} 
  where $N$ is the normal bundle of $\mathcal{Z}_\mu$ in the ambient space.   
     When we add this with  $(\ref{500})$ we get: 
    \begin{align}
   & (\pi\circ i)^*\pi_{*} i_{\mu*} (ev^*_{n+1}(F) \otimes \ozzr ) =  p_1^* (\pi_{*} i_{\mu*} (ev^*_{n+1}(F) \otimes \ozzr ) ) + \nonumber\\
  &   p_2^*(\pi_{*} i_{\mu*} (ev^*_{n+1}(F) \otimes \ozzr ))+  ev_{n+1}^* F \otimes(1-L_+L_- ) ,  \label{nodp8} 
    \end{align}
   as stated.

      \begin{tetacc} \label{tetaccc}
    The following hold:
   \begin{align}
   a.\quad & \pi^* \CC =   \mathcal{C}_{g,n+1,d}
    \cdot \prod_{j=1}^n\prod_{\delta=1}^{i_{\mu_j}} \mathcal{C}^{\mu_j}_\delta\left(-ev_{n+1}^*(F_{\delta\mu_j})\otimes \sigma_{j*}\mathcal{O}_{\mathcal{D}_j}\right)   \nonumber\\
  & \cdot \prod_\mu\prod_{\delta=1}^{i_\mu}\mathcal{C}^\mu_\delta\left(-ev_{n+1}^*(F_{\delta\mu})\otimes (i_{\mu*}\ozzr\right) .\label{nodalpull}\\
  b.\quad & p^* (i^{red}_\mu )^*\pi^* \CC =   \left( p_1^* \mathcal{C}^\mu_{g_1,n_1+1,d_1} \cdot p_2^* \mathcal{C}^\mu_{g_2,n_2+1, d_2} \right) \cdot \nonumber\\
     & \cdot (ev^*_+\times ev^*_-)\Delta_{\mu*}\left(\prod_{\delta=1}^{i_{\mu}} \mathcal{C}^\mu_\delta((q^*F_{\delta\mu})_\mu)\otimes (1- L_+L_-)))\right) .\label{restr3} \\
    c.\quad &p^* (i_\mu^{irr})^* \pi^* \CC = \nonumber\\
      &= \mathcal{C}^\mu_{g-1,n+2,d}\cdot (ev^*_+\times ev^*_-)\Delta_{\mu*}\left(\prod_{\delta=1}^{i_{\mu}} \mathcal{C}^\mu_\delta((q^*F_{\delta\mu})_\mu)\otimes (1- L_+L_-))\right) .
     \end{align}
     \end{tetacc}
  Proof: the equalities $(\ref{nodalpull})$ and $(\ref{restr3})$ are immediate consequences of $(\ref{nodp3})$ and $(\ref{nodp4})$ and of the multiplicativity of the classes $\CC$. 
 We will use $(\ref{nodalpull})$ in a different form, transforming the product into a sum:
  \begin{align}
   & \pi^* \CC = \mathcal{C}_{g,n+1,d}\cdot \prod_{j=1}^n \prod_{\delta=1}^{i_C}\left(1+ \mathcal{C}^{\mu_j}_{\delta}\left(- ev_{n+1}^*(F_{\delta\mu_j})\otimes \sigma_{j*}\oddi)\right) -1\right) \nonumber\\
  &\prod_\mu\left(1+ \prod_{\delta=1}^{i_{C_\mu}}\mathcal{C}^\mu_\delta\left(- ev^*_{n+1}(F_{\delta\mu})\otimes i_{\mu*}\ozzr\right)-1\right) =\nonumber\\
   & \mathcal{C}_{g,n+1,d} + \sum_{j} \mathcal{C}_{g,n+1,d}\cdot \prod_{\delta=1}^{i_{C_{\mu_j}}}\left(\mathcal{C}^{\mu_j}_\delta\left(- ev_{n+1}^*(F_{\delta\mu_j})\otimes \sigma_{j*}\oddi\right) -1\right)+\nonumber\\
    & + \sum_\mu \mathcal{C}_{g,n+1,d}\cdot  \left(\prod_{\delta=1}^{i_{C_\mu}} \mathcal{C}^\mu_\delta\left(-ev_{n+1}^*(F_{\delta\mu})\otimes i_{\mu*}\ozzr\right)-1\right) .\label{nodalsum}
   \end{align}
 This happens because the classes $\mathcal{C}^\mu_\delta(..) -1$ are supported on $\mathcal{D}_i$ and $\mathcal{Z}$ and $\mathcal{D}_i\cdot \mathcal{D}_j =\mathcal{D}_i\cdot \mathcal{Z}_\mu =0$ if $i\neq j$ (we'll use the same trick in $(\ref{reldifsheaf})$ below).

  We conclude the section by doing a short Grothendieck-Riemann-Roch computation which will turn out useful in the next section:
  \begin{nodalclass2} \label{aaa}
 Let $F\in K^0(\ix)$. Then
     \begin{align}
 ch  \left( \pi_*i_{\mu *}( ev_{n+1}^* F \otimes  \mathcal{O}_{\mathcal{Z}_\mu} ) \right) =  \pi_* i_{\mu *}\left( ch ( ev_{node}^* F^{(0)}_\mu ) \cdot  Td^\vee(-L_+\otimes L_- ) \right).     \label{chhh}
  \end{align}
    \end{nodalclass2}
  Proof: recall that $F_\mu^{(0)}$ is the invariant part of the restriction of $F $ to $\ix_\mu$ under the action of $g_\mu$ and that $r(\mu)$ is the   order of the distinguished node on $\mathcal{Z}_\mu$. 
   We'll simply write $r$ throughout the proof. 
  
  We apply Toen's theorem $\ref{tggrr}$ to the map $\pi$ to get:
     \begin{align*}
     ch  \left( \pi_*i_{\mu *}( ev_{n+1}^* F \otimes  \mathcal{O}_{\mathcal{Z}_\mu} ) \right) =  I\pi_*\left[\widetilde{ch}\left( i_{\mu *}( ev_{n+1}^* F \otimes \mathcal{O}_{\mathcal{Z}_\mu})\right) \cdot  \widetilde{Td}( T_\pi ) \right].
     \end{align*}
     We follow closely Section 7 of \cite{tseng}: there are three types of contributions of the inertia orbifold of the universal family mapping to the main stratum $\ix_{g,n,d}$. The main stratum $\ix_{g,n+1,d}$ contributes  
    $$ \pi_*\left[ch\left( i_{\mu *}( ev_{n+1}^* F \otimes \mathcal{O}_{\mathcal{Z}_\mu})\right) \cdot  Td( T_\pi ) \right]. $$
    
    Using the fact that the class $ev_{n+1}^* F \otimes \mathcal{O}_{\mathcal{Z}_\mu}$ is supported on $\mathcal{Z}_\mu$ and the formula (\ref{reldifsheaf}) for $Td(T_\pi) = Td^\vee (\Omega_\pi)$ one sees that the expression above equals 
    \begin{align} \label{stsum}
    \pi_* i_{\mu *}\left( ch ( ev_{n+1}^* F ) \cdot  Td^\vee(-L_+\otimes L_- ) \right).
    \end{align}
    
    The contributions from the divisors of marked points are trivially $0$ because these do not intersect $\mathcal{Z}_\mu$. In addition 
    there are contributions from $r-1$ copies of the locus $\mathcal{Z}_\mu$ itself. The invariant part of the relative tangent bundle of the map $\pi$
    is given by $-L^\vee_+\otimes L^\vee_-$ while the moving part is $L_+^\vee\oplus L_-^\vee$. Using the identification of each copy of $\mathcal{Z}_\mu$
    given in \cite{tseng} we get contributions of the form 
    \begin{align}
    (\pi\circ i_\mu)_*\left[ \widetilde{ch}\left(i_\mu^* i_{\mu *} ( ev_{n+1}^* F )\right) \cdot \frac{ Td^\vee(-L_+\otimes L_- ) }{(1-\zeta^l e^{\psi_+})(1-\zeta^{-l}e^{\psi_-})}\right].
    \end{align}
    for $\zeta^r=1$ and $1\leq l\leq r-1$. Notice that the Euler class of $\mathcal{Z}_\mu$ 
    in $\ix_{g,n+1,d}$ is $(1-L_+)(1-L_-)$. Hence, using also Lemma $7.3.6$ of \cite{tseng} the expresion above equals  
    $$(\pi\circ i_\mu)_*\left[ \widetilde{ch}(ev_{node}^* F)\cdot Td^\vee(-L_+\otimes L_- )\right] .$$
    
    Now notice that each of the $r$ contributions come with a weight of $1/r$ due to the fact that  a curve in $\mathcal{Z}_\mu$ has $\mathbb{Z}_r$ worth more automorphism than a nodal curve on the base and that the summation $\widetilde{ch}(F)$ on all $r$ components kills the non-invariant part of $F$ 
    under the action of $g_\mu$. 
      This finishes the proof of the statement. 
   \section{Proofs of Theorems}
Proof of Theorem $\ref{loop11}$: this is an easy consequence of Tseng's result and of the commutativity of the operators $\Delta_\alpha$. \vspace{12pt} \\  
  Proof of Theorem $\ref{dilat}$: remember that $\BB$ is a product of $i_B$  multiplicative characteristic classes.  We'll prove the statement using induction on $i_B$. The case $i_B=0$ is trivial. Assuming the statement holds for $i_B-1$, we'll prove the infinitesimal version of the proposition for $i_B$. Namely assume the twisting class $\mathcal{B}_{i_B}$ to be:
 \begin{align*}
  \mathcal{B}_{i_B}=exp\left(\sum_{l\geq 1} v_l ch_l \pi_*\left(f(L_{n+1}^{-1})-f(1))\right) \right).
  \end{align*}
  We  compute:
   \begin{align}
  & \frac{\partial \mathcal{D}_{A,B,C}}{\partial v_l} \mathcal{D}^{-1}_{A,B,C}=\nonumber\\
   = &\sum_{d,n}\frac{Q^d \hbar^{g-1}}{n!} \left\langle \prod_{i=1}^n \mathbf{t}(\ops_i) \cdot  ch_l \pi_*\left(f(L_{n+1}^{-1})-f(1)\right)\cdot  \Theta_{g,n,d} \right\rangle_{g,n,d}.
   \end{align}

  To compute $ ch_l \pi_*\left(f(L_{n+1}^{-1})-f(1)\right)$ above we apply Toen's GRR to the morphism $\pi$ to get:
    \begin{align}
   ch\left( \pi_*(f(L_{n+1}^{-1})-f(1))\right) = I\pi_*\left(\widetilde{ch}(f(L_{n+1}^{-1})-f(1)) Td^\vee(\Omega_\pi ) \right).
     \end{align}

   Notice that  $\widetilde{ch}=ch$ because the last marked point is not an orbifold point. We have:
     \begin{align}
     \widetilde{ch}(f(L_{n+1}^{-1})-f(1))= f(e^{-\psi_{n+1}})-f(1). \label{27}
     \end{align}
 In our situation there are three strata on the universal curve which get mapped to $\xii$:
      \begin{itemize}
     \item The total space $\xiii$.
     \item The locus of marked points $\di$.
      \item The nodal loci $\mathcal{Z}_\mu$ where $\mu\neq 0$, i.e. the node is an orbifold point.
      \end{itemize}
 But the expression on the RHS in $(\ref{27})$ above is a multiple of $\psi_{n+1}$ and $\psi_{n+1}$ vanishes on the locus of marked points $\mathcal{D}_j$  and on the locus of nodes $\mathcal{Z}$. Hence only the total space contributes to GRR. Exact sequences very similar with $(\ref{ses1}), (\ref{sepor})$ in Section $5$ allow us to write the sheaf of relative differentials (see also \cite{tseng}):
\begin{align}
\Omega_{\pi} =L_{n+1}-  \oplus_{j=1}^n (\sigma_{j})_{*} \mathcal{O}_{\di} - i_* L . \label{diff}
    \end{align}
  Keeping in mind that the bundle $L$ defined in Section $5$ has trivial Chern character we get :
    \begin{align}
      Td^\vee(\Omega_\pi) = Td^\vee(L_{n+1})\prod_{j=1}^n Td^\vee(-\sigma_{j*} \mathcal{O}_{\di})Td^\vee(-i_* \mathcal{O}_{\mathcal{Z}}).
       \end{align}
     We now use the fact that $\psi_{n+1}\cdot \mathcal{D}_j =\psi_{n+1}\cdot \mathcal{Z} = 0$ (recall that $Td^\vee(L_{n+1})-1$ is a multiple of $\psi_{n+1}$) to rewrite the product above as a sum:
   \begin{align}
Td^{\vee}(\Omega_{\pi})= & Td^\vee(L_{n+1}) + \sum_{j=1}^n (Td^\vee(-\sigma_{j*} \mathcal{O}_{\di}) -1) + \nonumber\\
 + & Td^\vee (-i_* \mathcal{O}_{\mathcal{Z}})-1. \label{reldifsheaf}
   \end{align}
 The last $n+1$ summands are classes supported on $\mathcal{D}_j$ and $\mathcal{Z}$, so they are killed by the presence of $\psi$ in $f(e^{-\psi_{n+1}})-f(1)$. After all these cancelations we see that:
\begin{align}
ch\left( \pi_*(f(L_{n+1}^{-1})-f(1))\right) = \pi_*\left((f(e^{-\psi_{n+1}})-f(1)) \cdot Td^\vee(L_{n+1})\right).\label{2121}
 \end{align}
  $(\ref{2121})$ is a linear combination of kappa classes $K_{aj}= \pi_*(ev_{n+1}^*\varphi_a \psi_{n+1}^{j+1})$.
 Now we pull the correlators back on the universal orbicurve. It is essential here that the corrections in the $\CC$  and $\ops_j$ classes  are also supported on $\mathcal{D}_j$ and $\mathcal{Z}$ (as we can see from  $(\ref{nodalsum})$ ) and the presence of $\psi_{n+1}$ kills them. Therefore (we denote by $[f]_l$ the homogeneous part of degree $l$ of $f$):
 \begin{align}
 &\mathcal{D}^{-1}_{\mathcal{A},\mathcal{B},\mathcal{C}}\frac{\partial \mathcal{D}_{\mathcal{A,B,C}}}{\partial v_l} = \sum_{d,n,g}\frac{Q^d \hbar^{g-1}}{n!} \int_{[\ix_{g,n+1,d}]} \prod_{i=1}^n \left( \sum_{k_i\geq 0} \left(ev_i^*(t_{k_i})\cdot \ops_i^{k_i} \right)\right) \cdot \Theta_{g,n+1,d}\cdot\nonumber\\
  &\cdot\left[(f(e^{-\psi_{n+1}})-f(1)) \cdot Td^\vee(L_{n+1})\right]_{l+1}\cdot   \prod_{\beta=1}^{i_B}\mathcal{B}_\beta\left( -\frac{f_\beta(L_{n+1}^{-1})-f_\beta (1)}{L_{n+1}-1}\right) \nonumber\\
    & -\int_{\ix_{0,3,0}}  \varphi_a\psi_3^{m+1}(\cdots )- \int_{\ix_{1,1,0}} \varphi_a\psi_1^{m+1} (\cdots). \label{32}
 \end{align}
   The correction terms occur because the spaces $\ix_{0,3,0}$ and $\ix_{1,1,0}$ are not universal families. Notice that the first correction is always $0$ for dimensional reasons, and the second is $\neq 0$ only for $m= deg (\varphi_a) = 0$ (again for dimensional reasons). If we denote this contribution by $K_{l, i_B}$ then the constant $K_\mathcal{B}$ in Theorem \ref{dilat} equals $\prod_{i, l} e^{K_{l, i}}$. This will not play any role further. 

 So the "new" twisting by the class $\mathcal{B}_{i_B}$ has the same effect as the translation:
   \begin{align*}
 \mathbf{t}_{\mathcal{B}}(z)= \mathbf{t}(z) + z - z \prod_{\gamma=1}^{i_B }\mathcal{B}_\beta\left( -\frac{f_\beta(\mathbf{L}_z^{-1})-f_\beta(1)}{\mathbf{L}_z-1}\right),
    \end{align*}
   because both potentials satisfy the same differential equation. To see this differentiate the potential $\mathcal{D}_{\mathcal{A,B}}(\mathbf{t}_{\mathcal{B}}(z))$ in $v_l$:
   \begin{align}
 &\frac{\partial \mathcal{D}_{\mathcal{A,B}}(\mathbf{t}_{B}(z))}{\partial v_l} \mathcal{D}^{-1}_{\mathcal{A,B}} =  \sum_{d,n,g}\frac{Q^d \hbar^{g-1}}{n!} \int_{[\ix_{g,n+1,d}]} \prod_{i=1}^n \left( \sum_{k_i\geq 0} \left(ev_i^*(t_{k_i})\cdot \ops_i^{k_i} \right)\right)  \cdot\nonumber\\
  & \cdot \psi_{n+1}ch_l\left(\frac{f(L_{n+1}^{-1})-f(1)}{L_{n+1}-1}\right)\cdot \Theta_{g,n+1,d}\cdot \prod_{\beta=1}^{i_B}\mathcal{B}_\beta\left( -\frac{f_\beta(L_{n+1}^{-1})-f_\beta (1)}{L_{n+1}-1}\right).\label{33}
   \end{align}
 But:
 \begin{align}
   \psi_{n+1}ch_l\left(\frac{f(L_{n+1}^{-1})-f(1)}{L_{n+1}-1} \right)= \psi_{n+1}\left[\frac{f(e^{-\psi_{n+1}})-f(1)}{e^{\psi}-1}\right]_l = \nonumber\\
    \left[\psi_{n+1}\frac{f(e^{-\psi_{n+1}})-f(1)}{e^{\psi}-1}\right]_{l+1}  = \left[(f(e^{-\psi_{n+1}})-f(1)) \cdot Td^\vee(L_{n+1})\right]_{l+1}\label{34}
 \end{align}
 because
  \begin{align*}
  Td^\vee(L_{n+1}) = \frac{\psi_{n+1}}{e^{\psi_{n+1}}-1}.
  \end{align*}
  Plugging $(\ref{34})$ in $(\ref{33})$ we see that $(\ref{33})$ and $(\ref{32})$ are of exactly the same form. The potentials also satisfy the same initial condition at $\mathbf{v}=0$ by  the induction hypothesis.\vspace{12pt}\\
    Proof of Theorem $\ref{pola11}$: we'll prove that
        \begin{align}
   \mathcal{D}_{\mathcal{A,B,C}} = exp\left(\frac{\hbar}{2}\sum_{a,b,\ga,\gb,\mu}A^\mu_{a,\alpha;b;\beta}\partial_{a}^{\alpha,\mu}\partial_b^{\beta,\mu^I} \right)\mathcal{D}_{\mathcal{A,B}}     \label{44}   
    \end{align}
where  $A^\mu_{a,\alpha;b,\beta}$ are the coefficients of the expansion:
 \begin{align}
      \sum_{a,b} A^\mu_{a,\alpha; b,\beta}\varphi_{\alpha,\mu} \ops_+^a\otimes \varphi_{\beta,\mu^I} \ops_-^b & = \frac{\Delta_{\mu*}\left(\prod_{\gd=1}^{i_\mu}\mathcal{C}^\mu_\gd\left( (q^*F_{\gd\mu})^{(0)}_\mu \otimes(1-\mathbb{L}_z) \right)-1\right)}{-\psi_+ -\psi_-}\in \nonumber\\
   & \in H^*(\ix_{\mu} ,\mathbb{Q})[\ops_+]\otimes H^*(\ix_{\mu^I},\mathbb{Q})[\ops_-] .  \label{polcoeff}
      \end{align}
Here $\psi_+ = c_1(L_+)$, $\psi_- =c_1(L_-)$    and $\Delta_\mu :\ix_\mu \rightarrow \ix_\mu \otimes \ix_{\mu^I}$ is the composition $(Id\times \iota)\circ \Delta$ . The map:
  \begin{align*}
     \Delta_{\mu_*}  :  H^*(\ix_\mu ,\mathbb{Q})\rightarrow H^*(\ix_{\mu} ,\mathbb{Q})\otimes H^*(\ix_{\mu^I},\mathbb{Q})
   \end{align*}  
   extends naturally to a map, which we abusively also call $\Delta_{\mu *}$ :
      \begin{align*}
     \Delta_{\mu_*}  :  H^*(\ix_\mu ,\mathbb{Q})[z] \rightarrow H^*(\ix_{\mu} ,\mathbb{Q})[\ops_+]\otimes H^*(\ix_{\mu^I},\mathbb{Q})[\ops_-],
     \end{align*} 
  by mapping $z\mapsto \psi_+\otimes 1 + 1\otimes\psi_-$   
  and the RHS of $(\ref{polcoeff})$ should be understood in this way. 
       
       We'll prove $(\ref{44})$ using induction on the total number $\sum_\mu i_\mu$ of twisting classes $\mathcal{C}_\delta^\mu$. If $\sum i_{\mu}=0$ then the equality is trivial. Let now $\sum i_{\mu}\geq 1$. Assuming $(\ref{44})$ to be true for $\sum i_{\mu}-1$, we'll prove the infinitesimal version of the theorem for $\sum i_{\mu}$. More precisely fix an $\mu_0$ and let the multiplicative class $\mathcal{C}^{\mu_0}$ (we omit the lower index) be of the form :
   \begin{align}
   \mathcal{C}^{\mu_0} ( E ) = exp\left( \sum_l w_l ch_l(E)\right) .\label{205}
   \end{align}
  As we vary the coefficients $w_l$ we obtain a family of elements in the Fock space. We prove $(\ref{44})$ by showing that both sides satisfy the same differential equations  with the same initial condition.   Notice that the induction hypothesis ensures that both sides of $(\ref{44})$ satisfy the same initial condition at $\mathbf{w}=0$. Moreover $\partial \mathcal{D}_{\mathcal{A,B}}/\partial w_l = 0$ so on the RHS only the coefficients $A^{\mu_0}_{a,\alpha;b,\beta}$ depend on $w_l$. So if denote the RHS by $\mathcal{G} $ and differentiate it  we get:
    \begin{align}
   \frac{\hbar}{2}\sum_{a,b}\frac{\partial A^{\mu_0}_{a,\alpha;b,\beta}}{\partial w_l}\partial_a^{\alpha,\mu_0} \partial_b^{\beta,\mu_0^I} \mathcal{G} = \frac{\partial}{\partial w_l}\mathcal{G}.\label{300}
    \end{align}
     To compute $\partial A^{\mu_0}_{a,\alpha;b,\beta}/\partial w_l  $   we differentiate  in $w_l$ relation $(\ref{polcoeff})$ to get:
   
    \begin{align}
  &\sum_{a,\alpha ;b,\beta} \frac{\partial A^{\mu_0}_{a,\alpha;b,\beta}}{\partial w_l} \varphi_{\ga,\mu_0} \ops_+^a \otimes \varphi_{\beta,\mu_0^I} \ops_-^b = \frac{-1}{\psi_+ +\psi_-}\cdot\nonumber\\ 
    & \cdot \Delta_{\mu_0 *}\left(ch_l \left((q^*F)^{(0)}_{\mu_0}(1-L_+L_-)\right)\prod_{\delta=1}^{i_{\mu_0}}\mathcal{C}^{\mu_0}_\delta\left( (q^*F)^{(0)}_{\mu_0}(1-L_+ L_-) \right)\right).\label{206}
    \end{align}
  But:  
   \begin{align}
  ch_l( (q^* F)^{(0)}_{\mu_0} (1-L_+L_-))= \left[ch (q^*F)^{(0)}_{\mu_0} (1-e^{\psi_+ +\psi_-})\right]_l  , \label{207}
   \end{align}
 hence
  \begin{align}
  &\sum_{a,\alpha ;b,\beta} \frac{\partial A^{\mu_0}_{a,\alpha;b,\beta}}{\partial w_l} \varphi_{\ga,\mu_0} \ops_+^a \otimes \varphi_{\beta,\mu_0^I} \ops_-^b = \frac{- 1}{\psi_+ +\psi_-}\cdot\nonumber\\ 
    & \cdot \Delta_{\mu_0 *} \left(\left[ch (q^*F)^{(0)}_{\mu_0} (1-e^{\psi_+ +\psi_-})\right]_l \prod_{\gamma=1}^{i_{\mu_0}}\mathcal{C}^{\mu_0}_\delta\left( (q^*F)^{(0)}_{\mu_0}(1-L_+ L_-) \right) \right).\label{240}
   \end{align}
 
  Below we prove that $\mathcal{D}_{\mathcal{A,B,C}}$ satisfies the same second order differential equation.
  The partial derivative of $ \mathcal{D}_{\mathcal{A,B,C}}$ with respect to $w_l$ equals:
       \begin{align}
      & \mathcal{D}^{-1}_{\mathcal{A,B,C}}\frac{\partial \mathcal{D}_{\mathcal{A,B,C}}}{\partial w_l} =\nonumber\\
   = &\sum_{d,n}\frac{Q^d \hbar^{g-1}}{n!}\left\langle\mathbf{t}(\ops_1 ),\ldots ,\mathbf{t}(\ops_n )  ;ch_l\pi_*(ev_{n+1}^*(F)\otimes i_{\mu_0*}\mathcal{O}_{\mathcal{Z}_\mu} )
\cdot\Theta_{g,n,d} \right\rangle_{g,n,d} .  \label{200}
       \end{align}
  Lemma $\ref{aaa}$ shows that:
  \begin{align}
ch_l\pi_*(ev_{n+1}^*(F)\otimes i_{\mu_0*}\ozz )  = \pi_* i_{\mu_0*}\left[  ev_{node}^* ch(F_{\mu_0}^{(0)}) \cdot \frac{e^{\psi_+ +\psi_-}-1}{\psi_+ +\psi_-}\right]_{l-1} .\label{202}
  \end{align}
Using $(\ref{202})$ and the formula :
  \begin{align*}
  \int_{[\ix_{g,n,d}]}(\pi_* i_* a) \cdot b = \int_{[\mathcal{Z}]} a \cdot(\pi\circ i)^* b
   \end{align*}
   we pullback the RHS of $(\ref{200})$ on $\mathcal{Z}$. Moreover we use Proposition $\ref{mainpllbc}$ to pullback the correlators on the factors $\ix_{g_1,n_1+1,d_1}\times \ix_{g_2,n_2+1,d_2}$.

       The classes $[\ix_{g,n,d}]^{tw}$ pullback as in formulae $(\ref{29}) , (\ref{30})$.
   As a consequence we see that if we define the  coefficients $A^{\mu_0 , l}_{a,\alpha;b,\beta}$  by:
    \begin{align}
    &\sum_{a,b,\alpha,\beta}A^{\mu_0,l}_{a,\alpha;b,\beta}\varphi_{\alpha,\mu_0} \ops_+^a \otimes \varphi_{\beta,\mu_0^I} \ops_-^b = \nonumber\\
&= \Delta_{\mu_0*}\left(\left[ch(q^*F)^{(0)}_{\mu_0}\cdot \frac{e^{\psi_+ +\psi_-}-1}{\psi_+ + \psi_-}\right]_{l-1}\left(\prod_{\delta=1}^{i_{\mu_0}}\mathcal{C}_\delta((q^*F)^{(0)}_{\mu_0}\otimes (1- L_+ L_-)\right)\right) ,\label{209}
      \end{align}
  we can express  $(\ref{200})$ as:
     \begin{align}
     & \mathcal{D}^{-1}_{\mathcal{A,B,C}} \frac{\partial \mathcal{D}_{\mathcal{A,B,C}}}{\partial w_l}= \sum_{g_i,n_i,d_i}\frac{Q^{d_1 +d_2}\hbar^{g_1+g_2-1}}{n_1 ! n_2 !}\cdot\nonumber\\
     & \cdot\sum_{a,b,\alpha,\beta}\frac{1}{2}\la \mathbf{t},\ldots ,\mathbf{t}, A^{\mu_0,l}_{a,\alpha;b,\beta}\varphi_{\alpha,\mu_0} \ops_+^a ;\Theta_{g_1,n_1+1,d_1}\ra_{g_1,n_1+1,d_1}  \times \nonumber\\
      &\times \la \mathbf{t},\ldots ,\mathbf{t},\varphi_{\beta,\mu_0^I} \ops_-^b ;\Theta_{g_2,n_2+1,d_2}\ra_{g_2,n_2+1,d_2} + \nonumber\\
       & + \sum_{\substack{g,n,d\\a,b,\ga,\gb}}\frac{1}{2}\frac{Q^{d}\hbar^{g-1}}{n!}\la \mathbf{t},\ldots ,\mathbf{t}, A^{\mu_0,l}_{a,\alpha;b,\beta}\varphi_{\alpha,\mu_0} \ops_+^a , \varphi_{\beta,\mu_0^I} \ops_-^b ; \Theta_{g-1,n+2,d}\ra_{g-1,n+2,d}.
     \end{align}

  Hence the generating function $\mathcal{D}_{\mathcal{A,B,C}}$ satisfies the equation:
    \begin{align}
    \frac{\partial\mathcal{D}_{\mathcal{A,B,C}}}{\partial w_l} = \frac{\hbar}{2}\sum_{a,b}  A^{\mu_0,l}_{a,\alpha;b,\beta}\partial_a^{\alpha,\mu_0}\partial_b^{\beta,\mu_0^I} \mathcal{D}_{\mathcal{A,B,C}}.\label{204}
    \end{align}
 Comparing  $(\ref{240})$ with $(\ref{209})$ we see that 
   \begin{align}
  \frac{\partial A^{\mu_0}_{a,\alpha;b,\beta}}{\partial w_l} =  A^{\mu_0,l}_{a,\alpha;b,\beta}.
    \end{align}
    
  Therefore both sides of  $(\ref{44})$ satisfy the same PDE. The theorem follows. 
\begin{coord1}
\emph{According to \cite{Tom} (pages $91-95$) this change of generating function corresponds to a change of polarization, namely we regard the potential $\mathcal{D}_{\mathcal{A,B,C}}$ as an element of the Fock space  $\mathcal{H}_\mathcal{C}=\mathcal{H}_+ \oplus \mathcal{H}_{-,\mathcal{C}}$ . The corresponding element in $\mathcal{H}=\mathcal{H}_+ \oplus \mathcal{H}_-$ with the usual polarization is $\mathcal{G}$. If $\{ q_a^{\alpha,\mu}, p_b^{\beta,\mu}\}$, $\{\overline{q}_a^{\alpha,\mu},\overline{p}_b^{\beta,\mu} \}$ are Darboux coordinates systems on $\mathcal{H}$, respectively $\mathcal{H}_{\mathcal{C}}$ then this change of polarization is given in coordinates by:}
\end{coord1}

    \begin{align}
   & p_b^{\beta,\mu} = \overline{p}_{b}^{\beta,\mu} ,\nonumber\\
   & \overline{q}_a^{\alpha,\mu} = q_a^{\alpha,\mu} - \sum_{a,b}A^\mu_{a,\alpha; b,\beta} p_b^{\beta,\mu} .\label{coordinates}
    \end{align}
  \begin{coord2} \label{totex}
   \emph{ Let $\ix$ be a manifold and let  $\mathcal{C}(\pi_* i_* \mathcal{O}_{\mathcal{Z}}) = Td( - \pi_* i_* \mathcal{O}_{\mathcal{Z}})^\vee$.  Then $A_{a,\alpha;b,\beta}$ don't depend on $\alpha$ or $\beta$ and we have:}
  
    \end{coord2}   
   
    \begin{align*}
    \mathcal{C}(1-L_+L_-) = Td^\vee(L_+L_-) = \frac{-\psi_+ - \psi_-}{1-e^{\psi_+ +\psi_-}}.
    \end{align*}
  This gives:  
    \begin{align*}
    \sum_{a,b} A_{a,\alpha, b,\beta}\psi^a \psi^b = \frac{1}{\psi_+ + \psi_-} - \frac{1}{e^{\psi_+ +\psi_-} -1}. 
    \end{align*}
   According to \cite{Tom} the expansion of :
    \begin{align*}
    \frac{1}{1-e^{\psi_+ + \psi_-}}= \sum_{k\geq 0} \frac{e^{k\psi_+}}{(1-e^{\psi_+})^{k+1}}(e^{\psi_-}-1)^k
     \end{align*} 
     gives a Darboux basis on $\mathcal{H}_\mathcal{C}$ in the sense of Theorem $\ref{pola11}$ i.e. $\varphi_a\frac{e^{k\psi_+}}{(1-e^{\psi_+})^{k+1}}$ span $\mathcal{H}_{-}$. 
     
\section{Quantum fake Hirzebruch-Riemann-Roch}
 As a first application we recover the quantum Hirzebruch-Riemann-Roch theorem of \cite{Tom}, which expresses the potential of the fake cobordism 
theory in terms of the cohomological one. Throughout this section, $X$ will be a compact complex manifold. 

 We first briefly review some basic background facts on complex-oriented cohomology theories. A more detailed review is given in \cite{Tom}.   
\begin{defcoh}
{\em A complex-oriented cohomology theory is a multiplicative cohomology theory $E^*$ together with a choice of element $u_E\in E^2(\mathbb{CP}^\infty)$ such that if $j:\mathbb{CP}^1\to \mathbb{CP}^\infty$ is the inclusion then $j^*(u_E)$ is the standard generator of $E^2(\mathbb{CP}^1)$.}
 \end{defcoh}
 We denote the ground ring by $R_E:=E^*(pt.)$. One can define Chern classes satisfying the usual axioms such that $j^*(u_E)$ is the first Chern class of the Hopf bundle. 
 The Chern-Dold character is the unique multiplicative natural transformation:
   \begin{align*}
  ch_E : E^*(X)\to H^*(X,R_E)   
     \end{align*}
   which is the identity if $X=\{pt\}$.
   
 In particular $ch_E(u_E)$ is a power series in $z$, where $z$ is the standard orientation of $H^*(X,R_E)$. We denote it $u_E(z)$. The Todd class is the unique multiplicative class which for a line bundle $L$ is:
 \begin{align*}
Td_E(L) : = \frac{c_1(L)}{u_E(c_1(L))} . 
  \end{align*}
   
  We now fix the cohomology theory to be complex cobordism $MU^*$. For a given $i$, $MU^i(X)$ is defined as:
 \begin{align*} 
MU^i(X):= lim_{j \to \infty}[\Sigma^j X, MU(i+j)],
  \end{align*}
 where $[$ , $]$ denotes homotopy classes of maps, $\Sigma^j X$ is the iterated reduced suspension of $X$ and $MU(k)$ are the Thom spaces. 
 
 Cobordism is universal among complex-oriented cohomology theories in the following sense: for any other cohomology $(E, u_E)$ there is a unique natural transformation $MU\to E$ which maps $u$ to $u_E$ (we will write $u,R$ etc. instead of $u_{MU},R_{MU}$). If $X$ has complex dimension $n$, $MU^i(X)$ can be identified with the complex bordism group $MU_{2n-i}(X)$. This is Poincar\'e duality for complex cobordism and bordism. The image of $u$ under the Chern-Dold map is a formal power series $u(z)$, where $z$ is the first Chern class of the universal line bundle.
 
The ground-ring of the cobordism is $R:=MU^*(pt.)=\mathbb{C}[p_1,p_2,\ldots ]$ (we tensored with $\mathbb{C}$) where $p_i$ is the class of the map $\mathbb{CP}^{i}\to pt.$. 
For a l.c.i. map $f:X\to Y$ there is a push-forward $f_*$ and a Hirzebruch-Riemann-Roch theorem which says the diagram:
   \begin{displaymath}
     \begin{CD}
    MU^*(X)  @>ch_{MU}\cdot Td(T_f) >> H^*(X, R)\\
      @V f_* VV                 @V  f_* VV                          \\
 MU^*(Y) @> ch_{MU}>>  H^*(Y,R)
             \end{CD}
 \end{displaymath}
 is commutative. We define ``fake'' cobordism-valued Gromov-Witten invariants to be given by the above theorem applied to the morphisms $X_{g,n,d}\to \{pt\}$.
  
 Denote by $\mathcal{T}_{g,n,d}$ the virtual tangent bundle to $X_{g,n,d}$. The genus-$g$ descendant cobordism-valued potential (called ``extraordinary potential'' in \cite{Tom}) is defined as:
 
 \begin{align*}
\mathcal{F}_{MU}^g : = \sum_{d,n}\frac{Q^d}{n!}\int_{[X_{g,n,d}]}\prod^n_{i=1}\left(\sum_{k\geq 0}ch_{MU}(ev_i^* t_k u(\psi_i)^k )\right)\cdot Td_{MU}(\mathcal{T}_{g,n,d}). 
 \end{align*}
 It is a formal function of
  \begin{align*}
 \mathbf{t}(u):=\sum_{k\geq 0} t_k u^k \in MU^*(X)[[u]]
   \end{align*}
 which takes values in the ring $R[[Q]]$. The total extraordinary potential is
  \begin{align*}
 \mathcal{D}_{MU}:= exp\left(\sum^\infty_{g=0} \hbar^{g-1}\mathcal{F}_{MU}^g\right)  .
  \end{align*}  
  
We define $\mathcal{U}$ to be the space:
 \begin{align*}
 \mathcal{U}:= MU^*(X,\mathbb{C}[[Q]])[[u]].
  \end{align*}
  The symplectic form on $\mathcal{U}$ is:
  \begin{align*}
  \Omega_{MU}\left(\mathbf{f},\mathbf{g} \right):= \oint_{z=0} \left(\mathbf{f}(u(z),\mathbf{g}(u(-z)\right)_{MU}dz 
   \end{align*}
  with the pairing : 
 \begin{align*}
\left(\ga, \gb \right)_{MU} =\int_X ch_{MU}(\ga)\cdot ch_{MU}(\gb)\cdot Td_{MU}(T_X) .
 \end{align*} 
 The space $\mathcal{U}_+$  of the polarization on $\mathcal{U}$ is defined to include all power series in $u$. If we expand:
  \begin{align*}
 \frac{1}{u(-x-y)} = \sum_{k\geq 0} u^k(x) v_k(u(y)) 
   \end{align*}
  then $\mathcal{U}_-$ is defined as the span of all $\phi_\alpha v_k(u)$ for all $k\geq 0$, $\phi_\alpha\in MU^*(X)$. It is shown in \cite{Tom} that these two subspaces realize a polarization of $\mathcal{U}$.
To show how the extraordinary potential is related to the cohomological one we define a modification of $\mathcal{H}$:
 \begin{align*}
\mathcal{H}_{MU}: = H^*(X,R[[Q]])((z)).
 \end{align*}   
 The pairing and symplectic form on $\mathcal{H}_{MU}$ (henceforth denoted $\mathcal{H}$) are defined in the obvious way.
  The map:
  \begin{align*}
 &\widetilde{ch}_{MU}: \mathcal{U} \to \mathcal{H} ,\\
 &\sum_k t_k u^k \mapsto \sqrt{Td_{MU}(T_X)} \left(\sum_k ch_{MU}(t_k) u^k(z)\right)
  \end{align*}
 is a symplectomorphism which maps $\mathcal{U}_+$ to $\mathcal{H}_{+}$, but it doesn't map $\mathcal{U}_-$ to $\mathcal{H}_-$.
Let 
 \begin{align}
\mathbf{q}(z) = \sqrt{Td_{MU}(T_X)}( \mathbf{t}(z) + u(-z)). \label{3434}
 \end{align} 
We regard $\mathcal{F}^0_{MU}, \mathcal{D}_{MU}$ as functions of $\mathbf{q}(z)$ (hence a function on $\mathcal{H}_+$) via the identifications above. Let 
$\widehat{\nabla}$ be the quantized linear symplectic transformation 
 \begin{align*}
 \widehat{\nabla} := exp\left( A_{a,\ga;b,\gb}g^{\ga\gb}\partial_a^\ga\partial_b^\gb\right)
  \end{align*}
  with $A_{a,\ga;b,\gb}$ given in Example 4.2. Let:
\begin{align}
 \Delta:= exp\left(\sum_{m\geq 0}\sum^{dim(X)}_{l=0} s_{2m-1+l}\frac{B_{2m}}{(2m)!}ch_l(T_X)z^{2m-1} \right) \label{3436}
 \end{align}
 where $s_k$ are defined by:
  \begin{align*}
 exp\left(\sum_{k\geq 1}s_k\frac{x^k}{k!} \right) = \frac{x}{u(x)} \in H^*(X, R). 
   \end{align*}
 \begin{indexx0} \label{3500}
We then have:
 \begin{align*}
\mathcal{D}_{MU} \approx \widehat{\nabla} \widehat\Delta \mathcal{D} .
  \end{align*}
 \end{indexx0}
  The proof will be a consequence of the description of the virtual tangent bundles to $X_{g,n,d}$ as linear combinations of classes of type
$\mathcal{A,B,C}$:
\begin{tangvirtual}\label{tggv}
 \begin{align}
 \mathcal{T}_{g,n,d}:=\pi_*ev_{n+1}^*(T_X - 1)- \pi_*(L^{-1}_{n+1}-1) - \left(\pi_*i_*\mathcal{O}_\mathcal{Z} \right)^\vee .\label{tgg}
  \end{align}
  \end{tangvirtual}
  Proof: we follow closely the computation in the dissertation \cite{Tom}. However the proof there, while leading to the same formula, is a bit imprecise in assuming that $L_{n+1}$ restricted to $\mathcal{Z}$ is the trivial line bundle. Recall that $\mathcal{Z}$ is the nodal locus in the universal family and that it is parametrized by $\widetilde{\mathcal{Z}}$ which is a fiber product of moduli spaces of maps of lower genus. The gluing map $\widetilde{\mathcal{Z}}\to \mathcal {Z}$ is generically $2$ to $1$.  The symmetry on $\widetilde{\mathcal{Z}}$ permuting the two marked points which become the node after gluing acts non-trivially on the fibers of $i^* L_{n+1}$ above the fixed point locus. Hence $i^*L_{n+1}$ is a non-trivial (orbi)bundle on $\mathcal{Z}$. We denote it by $L$. More precisely  let $L'$ be the $\mathbf{Z}_2 $ equivariant line bundle on $\widetilde{\mathcal{Z}}$ which is $\widetilde{\mathcal{Z}}\times \mathbb{C}$ as a set and on which $-1\in \mathbb{Z}_2$ acts by :
  \begin{align*}
 (x, v)\mapsto (-1\cdot x , -v )\quad for \quad x\in \widetilde{\mathcal{Z}} ;\quad v\in \mathbb{C} ;
   \end{align*}
 Then $L= L'/\mathbf{Z}_2$.  
 
Let $(C, x_1,\ldots x_n)$ be a point in $X_{g,n,d}$ and let $D$ be the divisor of marked points $D =\mathcal{D}_1+\ldots + \mathcal{D}_n$. Then (see $\cite{Tom}$ and the references therein):
\begin{align}
 \mathcal{T}_{g,n,d} = &\pi_*(ev^*T_X) - H^0(C, \Omega^\vee_\pi(-D)) + H^1(C, \Omega^\vee_\pi(-D))\nonumber\\
  = &\pi_*(ev^*T_X) - \pi_*(\Omega_\pi^\vee (-D) ).\label{tg01}
    \end{align}
   Roughly the first summand acounts for deformations of the map, the second for infinitesimal automorphisms of the curve $(C, x_1, \ldots , x_n)$ and the third for deformations of the complex structure of $C$ and smoothing of the nodes.  
Denote by $\omega_\pi$ the dualizing sheaf of the universal family. According to \cite{Tom} we have the exact sequence:
 \begin{align}
 0\to \omega_\pi \to L_{n+1} \to \oplus_j\sigma_{j*}(\mathcal{O}_{\mathcal{D}_{j}})\to 0 .\label{sepor}
 \end{align}
 Using Serre duality and the relation given by the above exact sequence the second summand in $(\ref{tg01})$ becomes: 
\begin{align}
- \pi_*(\Omega_\pi^\vee (-D) ) = \left[\pi_*(\Omega_\pi (D)\otimes \omega_\pi)\right]^\vee= \left[\pi_*(\Omega_\pi\otimes L_{n+1})\right]^\vee.  \label{tg02}
 \end{align}
 There is an exact sequence:
\begin{align}
0\to \Omega_\pi \to \omega_\pi \to i_* L \to 0. \label{ses1}
\end{align}

First notice that $\Omega_\pi$ and $\omega_\pi$ coincide away from $\mathcal{Z}$. Near a point of $\mathcal{Z}$ the map $\pi$ can be described locally by:
 \begin{align*}
 \pi: (z,x,y) \to (z,xy)
 \end{align*}
where $z$ is a (vector) coordinate on $\widetilde{\mathcal{Z}}$ viewed as an orbifold chart for $\mathcal{Z}$ and the symmetry $-1\in \mathbb{Z}_2$ interchanges $x$ and $y$.  Locally sections of $\omega_\pi$ have the form 
\begin{align*}
 f(z,x,y)\frac{dx\wedge dy}{d(xy)}
\end{align*}  
and of $\Omega_\pi$ are of the form $g(z,x,y)dx + h(z,x,y)dy$ where we impose the relation $xdy+ydx=0$. There is a natural inclusion:
\begin{align*}
& \Omega_\pi \to \omega_\pi  \\
& g(z,x,y)dx + h(z,x,y)dy \mapsto \left(xg(z,x,y) -yh(z,x,y)\right)       \frac{dx\wedge dy}{d(xy)}.
\end{align*}

Sections in the cokernel are represented by elements of the form :
\begin{align*}
\alpha(z)\frac{dx\wedge dy}{d(xy)}.
\end{align*}
This is identified with $i_* L $ because the symmetry acts non-trivially on  $dx\wedge dy$. This establishes $(\ref{ses1})$. 

We now use  $(\ref{ses1})$ to rewrite $\Omega_\pi = \omega_\pi - i_*L$ and then plug in $(\ref{tg02})$:  
\begin{align}
\left[\pi_*(\Omega_\pi\otimes L_{n+1})\right]^\vee = \left[\pi_*(\omega_\pi\otimes L_{n+1})\right]^\vee - \left[\pi_*(i_*(L)\otimes L_{n+1})\right]^\vee .\label{tg03}
    \end{align}
The first term in $(\ref{tg03})$ equals $-\pi_*[L^{-1}_{n+1}]$ by Serre duality again. Replacing in $(\ref{tg03})$ we get:   
          \begin{align}
\left[\pi_*(\Omega_\pi\otimes L_{n+1})\right]^\vee = -\pi_*[L^{-1}_{n+1}]-\left[\pi_*i_*(L\otimes i^*L_{n+1})\right]^\vee \label{tg04} 
    \end{align}                
   But $i^*L_{n+1}= L$ and $L^2=1$. Hence the last term in $(\ref{tg04})$ is $- \left(\pi_*i_*\mathcal{O}_\mathcal{Z} \right)^\vee$. Formula $(\ref{tgg})$ then follows by plugging $(\ref{tg04})$ in $(\ref{tg01})$.\\
 
  Proof of Theorem $\ref{3500}$: we regard the Todd class $Td_{MU}$ as a family of multiplicative classes depending on the parameters $s_i$. Then the twisting theorems apply:
\begin{itemize}
\item twisting by $Td_{MU}\left(\pi_*ev_{n+1}^*(T_X - 1)\right)$ corresponds to acting by the operator $\widehat\Delta$ on the potential $\mathcal{D}$ according to Remark $\ref{rem333}$;
\item twisting by $Td_{MU}\left(- \pi_*(L^{-1}_{n+1}-1) \right)$ accounts for the dilaton shift $(\ref{3436})$ according to Theorem $\ref{dilat}$;
\item twisting by the class $Td_{MU}\left(- \left(\pi_*i_*\mathcal{O}_\mathcal{Z} \right)^\vee\right)$: according to the proof of Theorem $\ref{pola11}$ and Example $\ref{totex}$ is tantamount to acting on the potential by the operator $\widehat{\nabla}$ . 
 \end{itemize}
  
By looking only at genus $0$  we easily deduce the following:
  \begin{indexx}
The graph of the generating series $\mathcal{F}^0_{MU}$, viewed as a formal function of $\mathbf{q}(z)$ with respect to the polarization
 \begin{align}
\mathcal{H}_{MU}=\mathcal{H}_+ \oplus \{\phi_\alpha v_k(u(z))\vert k\geq 0, \phi_\alpha\in H^*(X,R)\} 
 \end{align} 
 is a Lagrangian cone $\mathcal{L}_{MU}$. It is obtained from the cohomological cone $\mathcal{L}^H$ after rotating by the symplectic transformation $\Delta$.
    \end{indexx} 
  \section{Applications to the GW theory of $X\times B\mathbb{Z}_m$}   
 In this section we apply the results to the Gromov-Witten theory of the orbifold $X\times B\mathbb{Z}_m$, where $X$ is a smooth complex manifold. The motivation lies in the study of the quantum K-theory of $X$. The results in this section are used in Section 8 in \cite{gito}. 
 
   Let $G$ be a finite group which acts trivially on $X$ and let $\ix=X\times BG$, the stack theoretic quotient.   
  We denote by $[\gamma_i]$ the conjugacy class of $\gamma_i\in G$ and by $C(\gamma)$ the centralizer of $\gamma$. The inertia stack of $X/G$ is the disjoint union $\coprod_{i} ([\gamma_i], X/C(\gamma_i))$. Therefore :
    \begin{align*}
    H^*(I(X/G),\mathbb{C}) = \oplus_{[\gamma_i]} H^*(X,\mathbb{C}).
     \end{align*} 
   Denote by $e_{[\gamma_i]} : = 1 \in H^*(([\gamma_i], pt/C([\gamma_i])))$. A basis of $H^*(([\gamma_i], X/C(\gamma_i)))$ is given by $\varphi_a\times e_{[\gamma_i]}$, where $\{\varphi_a\}$ is a basis of $H^*(X,\mathbb{C})$. The Poincar\'e pairing is given by : 
    \begin{align*}
    ( \varphi_a \times e_{[\gamma_i]},\varphi_b\times e_{[\gamma_j]}) = \frac{\delta_{[\gamma_i][\gamma^{-1}_j]}}{\vert C(\gamma_i)\vert} \int_X \varphi_a \smile \varphi_b .
    \end{align*}
   The $J$ function is defined as:
   \begin{align}
   J_\ix (t, -z) = -z + t(z)+\sum_{n,d}\frac{Q^d}{n!}\phi_a\langle \frac{\widetilde{\phi^a}}{-z-\ops_1}, t(\ops_2), \ldots ,t(\ops_n)\rangle^{X/G}_{n,d}.
   \end{align}  
  where $\{\phi_a\}, \{\widetilde{\phi^a} \}$ are dual basis. We use results of \cite{jaki} to express the correlators in terms of correlators on $X_{0,n,d}$.
   In fact there is a finite degree map:  $(X\times BG)_{0,n,d,([\gamma_1],\ldots ,[\gamma_n])}\to X_{0,n,d}$. In \cite{jaki} it is shown the degree equals
       \begin{align*}
        \frac{\vert\chi^G_0(\mathbf{\gamma})\vert}{\vert G\vert} ,
      \end{align*}
   where
      \begin{align*}
  \chi_0^G(\mathbf{\gamma}): = \{(\sigma_1,\ldots ,\sigma_n) \vert 1=\prod^n_{j=1} \sigma_j, \sigma_j\in [\gamma_j] \quad for \quad all \quad j\}.    
      \end{align*}
      
      Since the $\ops$ classes in the correlators are pullbacks of $\psi$ classes from the coarse curve it follows that:
    \begin{align}
    \langle \prod_i \ops^{k_i}_i(ev_i^*(t_i \times e_{[\gamma_i]})\rangle^{X/G}_{0,n,d} =  \frac{\vert\chi^G_0(\mathbb{\gamma})\vert}{\vert G\vert} \langle \prod_i \psi^{k_i}_i ev_i^*(t_i) \rangle^{X}_{0,n,d} \label{4000}
    \end{align} 
 where $t_i\in H^*(X)$. \vspace{12pt}\\
  From now on, let $G=\mathbb{Z}_m$ and $\zeta$ a primitive $m$-th root of unity. Denote by $td_\zeta$ the multiplicative class defined for line bundles $L$ by:
      \begin{align*}
     td_\zeta (L) : = \frac{1}{1-\zeta e^{-c_1(L)}}.
     \end{align*}  
    We twist the cohomological potential of $\ix$ with $3$ types of twisting classes as follows: 
     \begin{itemize}
    \item   the type $\mathcal{A}$ classes we take to be:
                 \begin{align*}
             td(\pi_*ev^*(T_X)) \prod_{k=1}^{m-1} td_{\zeta^k}(\pi_*ev^*(T_X\otimes \mathbb{C}_{\zeta^k})).    
                  \end{align*}
     For a function $s(x)$, the Euler-Maclaurin asymptotics of $\prod^\infty_{r=1} e^{s(x-rz)}$ is  given by:
   \begin{align*}
  \sum^{\infty}_{r=1} s(x-rz) = \left(\sum^\infty_{r=1} e^{-rz\partial_x}\right) s(x) = \frac{z\partial_x}{e^{z\partial_x}-1}(z\partial_x)^{-1}s(x) \\
   = \frac{s^{(-1)}(x)}{z} -\frac{s(x)}{2} + \sum_{k=1}^\infty \frac{B_{2k}}{(2k)!}s^{(2k-1)}(x) z^{2k-1} , 
   \end{align*}
  where $s^{k} = d^ks/dx^k$, $s^{-1}$ is the antiderivative $\int_0^x s(t) dt$ and $B_{2k}$ are Bernoulli numbers. The effect of the type $\mathcal{A}$ twisting  is: 
     \begin{typeaa}  \label{crraa}
    The cone rotates by the loop group element: 
    \end{typeaa}  
      \begin{align*}
      \mathcal{L}^{tw}= \prod_{j=0}^{m-1}(\Box_{j}) \mathcal{L}_\ix ,
      \end{align*}
    where we think of $\mathcal{L}_\ix $ as a product of $m$ copies of $\mathcal{L}_X$ and each operator $\Box_j$ acts on the copy corresponding to the sector labeled by $g^j$. Let $[kj/m]$ denote the greatest integer less than $kj/m$. The operators in the statement are Euler-MacLaurin expansions of the products:
      \begin{align*}
    & \Box_0 = \prod_{i} \prod_{r=1}^{\infty}\frac{x_i-rz}{1-e^{-mx_i+mrz}}, \\
    & \Box_{j} =  \prod_{k=0}^{m-1}\prod_{i} \prod_{r=1}^{\infty}\frac{x_i-rz}{1-\zeta^k e^{-x_i+rz- (kj/m-[kj/m])z}}.
       \end{align*}    
  Proof: the corollary follows by 
 application of the main theorem of \cite{tseng} to the twisting data described above.
     \item  the type $\mathcal{B}$ classes :  
           \begin{align*}
      td(\pi_*(1-L^{-1}_{n+1})) \prod_{k=1}^{m-1} td_{\zeta^k}(\pi_*((1-L^{-1}_{n+1}\otimes ev^* \mathbb{C}_{\zeta^k})).       
            \end{align*}
      \begin{typebb}\label{typebb22}
    The dilaton shift changes from $q(z)= \mathbf{t}(z)-z$ to $q(z) = \mathbf{t}(z) +(1-e^{mz})$.  
      \end{typebb}  
   Proof: we apply Theorem $\ref{dilat}$ to the potential $\mathcal{F}$. \vspace{12pt}\\
    In our case $f_\beta = - ev_{n+1}^*(\mathbb{C}_\zeta )  \otimes L_{n+1}^{-1}$  we have:
   \begin{align*}
  \frac{f_\beta(L_{n+3}^{-1})-f_\beta(1)}{L_{n+3}-1} = \mathbb{C}_\zeta L_{n+3}^{-1}.
   \end{align*}
    So according to Theorem $\ref{dilat}$ (fix $\zeta$ to be primitive $m$-th root of unity) the translation is:
    \begin{align}
    \mathbf{t}(z)& : = \mathbf{t}(z) +z - z \prod_{k=0}^{m-1} Td_{\zeta^k} (-\mathbb{C}_{\zeta^k} \mathbf{L}_z^{-1}) =\nonumber\\
                 & : =  \mathbf{t}(z) +z - z \frac{1-e^{z}}{z }\prod_{k=1}^{m-1}(1-\zeta^k e^z) = \mathbf{t}(z) +z - (1-e^{mz}).
    \end{align}       
      \item  denote by $i_j$ the inclusion of the nodal locus $\mathcal{Z}_{g^j}$ for all $0\leq j\leq m-1$; the type $\mathcal{C}$ twisting we take to be
       \begin{align*}
     \prod_{j=0}^{m-1} \left[td^\vee(-\pi_*(i_{j*}\mathcal{O}_{\mathcal{Z}_{g^j}})) \prod_{k=1}^{m-1} td^\vee_{\zeta^k}(-\pi_*(i_{j*}\mathcal{O}_{\mathcal{Z}_{g^j}}\otimes ev^*\mathbb{C}_{\zeta^k}))\right].       
            \end{align*}               

      \begin{typecc} \label{typecc22}
   The nodal twisting changes the polarization as follows: in the sector $(\ix,g^j)$ of $I\ix$ the new Darboux basis is given by expansion of   
       \end{typecc}
      \begin{align*}
   \frac{1}{1-e^{d\psi_+ +d\psi_-}},
   \end{align*}
where $d = g.c.d. (j, m)$.   
      \end{itemize}   
    
Proof: according to Theorem $\ref{pola11}$, only the classes $td^\vee_{\zeta^k}$  for which $\mathbb{C}_{\zeta^k}$ is a trivial representation of 
$g^j$ give  nontrivial contributions to the twisting and there are $d$ worth of these.  
  The coefficients $A_{a,\alpha,b,\beta}^j$ are given by:
   \begin{align*}
 -\frac{\prod_{\zeta^{kj}=1} td_{\zeta^k}(1-L_+L_-) - 1}{\psi_+ +\psi_-} & = -\frac{1}{\psi_+ +\psi_-}\left(\frac{\psi_+ +\psi_-}{\prod_{\zeta^{kj}=1}(1-\zeta^k e^{\psi_+ +\psi_-})}-1\right)\\
  & =  \frac{1}{\psi_+ + \psi_-} - \frac{1}{e^{d\psi_+ +d\psi_-} -1}. 
   \end{align*}
   
   Then (see Example $\ref{totex}$ and \cite{Tom}) the Darboux basis is given by the expansion of $\frac{1}{1- e^{d\psi_+ +d\psi_-} } $.

\appendix
    \section{Grothendieck-Riemann-Roch for stacks}
      The main tool for proving Theorems $\ref{loop11}$, $\ref{dilat}$ and $\ref{pola11}$ is a generalization of Grothendieck-Riemann Roch theorem for morphisms of stacks due to B.To\"en (\cite{toen}). Before stating it we will introduce more notation:
  \begin{trace}
  \emph{Define $Tr:K^0(\ix )\to K^0(I\ix )$ to be the map:}
        \begin{align*}
      F\mapsto \oplus \lambda_i (g) F_i 
       \end{align*}
   \em{ on each component $(g,\ix_\mu )$ of the inertia stack, where $F_i$ is the decomposition of the $g$ action and $\lambda_i(g)$ is the eigenvalue of $g$ on $F_i$.}    
         \end{trace}
           \begin{cher}
  \emph {Define $\widetilde{ch}:K^0(\ix)\to H^*(I\ix )$ to be the map $ch \circ Tr$}.
         \end{cher}
      Now each vector bundle $E$ on $\ix $ restricts on each connected component $(g,\ix_\mu)$ of the inertia stack as the direct sum $E_{inv}\oplus E_{mov}$.
   \begin{todd}
     Define $\widetilde{Td}(E) :K^0(\ix )\to H^*(I\ix )$ to be the class:
    \end{todd}
       \begin{align*}
       \widetilde{Td}:= \frac{Td(E_{inv})}{ch(Tr\circ \lambda_{-1}(E_{mov})^\vee )}
       \end{align*}
    where $\lambda_{-1}$ is the operation in K-theory defined as $\lambda_{-1}(V):=\sum_{a\geq 0} (-1)^a \Lambda^a V$. In the following  theorem we assume the morphism $f$ factors as the composition of a smooth regular immersion followed by a smooth morphism. Then one can define $T_f$ as in the case of l.c.i. morphisms of manifolds. 
    \begin{Grr} \label{tggrr}
       Let $f:\ix\to \mathcal{Y}$ be a proper morphism of smooth Deligne Mumford stacks (over $\mathbb{C}$) with quasi-projective coarse moduli spaces. This induces a morphism $If:I\ix\to I\mathcal{Y}$. If $f$ factors as stated above we have:
       \begin{align}
   \widetilde{ch}(f_*E) = If_*\left(\widetilde{ch}(E)\widetilde{Td}(T_f) \right).
    \end{align}
\end{Grr}
   Restricting to the identity component $\mathcal{Y}$ of $I\mathcal{Y}$ we get:
    \begin{align} \label{tgrr}
ch (f_* E) = If_*\left(\widetilde{ch}(E)\widetilde{Td}(T_f)_{\vert If^{-1}\mathcal{Y}} \right).
    \end{align}
  The universal curve $\pi$, to which we apply Theorem $\ref{tggrr}$ is not necessarily a local complete intersection, so following \cite{tseng} we proceed as follows. The construction in \cite{abgot} provides a family of orbicurves 
  \begin{align}
  \widetilde{\pi}:\mathcal{U}\to \mathcal{M}
  \end{align}
  and an embedding $\ix_{g,n,d}\to \mathcal{M}$ satisfying the following properties:
   \begin{itemize}
  \item{The family $\mathcal{U}\to \mathcal{M}$ pulls back to the universal family over $\ix_{g,n,d}$.}
   \item{A vector bundle of the form $ev_{n+1}^*(E)$ extends to a vector bundle over $\mathcal{U}$.}
   \item{The Kodaira-Spencer map $T_m\mathcal{M}\to Ext^1(\mathcal{O}_{\mathcal{U}_m},\mathcal{O}_{\mathcal{U}_m})$ is surjective for all $m\in \mathcal{M}$.}
   \item{The locus $\mathcal{Z}\subset \mathcal{U}$ of the nodes of $\widetilde{\pi}$ is smooth and $\widetilde{\pi}(\mathcal{Z})$ is a divisor with normal crossings.}
   \item{The pull-back of the normal bundle $N_{\mathcal{Z}/\mathcal{U}}$ to the double cover $\widetilde{Z}$ given by choice of marked points at the node is isomorphic to the direct sum of the cotangent line bundles at the two marked points.} 
   \end{itemize}
 
 So technically we apply Grothendieck-Riemann Roch to $\widetilde{\pi}$ and then cap with the virtual fundamental classes $[\ix_{g,n,d}]^{tw}$. Therefore in the computations we assume the universal family $\pi$ satisfies the above properties.

Valentin Tonita, Kavli IPMU, University of Tokyo, 5-1-5 Kashiwa-no-Ha, Kashiwa City, Chiba 277-8583, Japan \\
\emph{E-mail:}\url{valentin.tonita@ipmu.jp}\\
Phone:0471366522
  \end{document}